\documentclass[reqno]{amsart}
\usepackage{mathtools,amssymb,amsthm,mathrsfs,graphicx}
\usepackage{empheq}

\usepackage{geometry}
\usepackage[numbers,sort&compress]{natbib}

\usepackage{todonotes}
\usepackage[draft=false,hidelinks]{hyperref}
\hypersetup{
  colorlinks   = true, 
  urlcolor     = blue, 
  linkcolor    = red, 
  citecolor   = blue 
}
\usepackage{esint}

\numberwithin{equation}{section}

\theoremstyle{plain}
\newtheorem{thm}{Theorem}[section]
\newtheorem{lem}[thm]{Lemma}
\newtheorem*{lem*}{Lemma}

\newtheorem{proposition}[thm]{Proposition}

\theoremstyle{definition}
\newtheorem*{defn}{Definition}
\newtheorem{assumption}[thm]{Assumption}

\theoremstyle{remark}
\newtheorem{remark}[thm]{Remark}

\allowdisplaybreaks

\newcommand{\R}{\mathbb{R}}
\newcommand{\C}{\mathbb{C}}

\begin{document}

\title[Parabolic systems]{Singular-degenerate parabolic systems with the conormal boundary condition on the upper half space}

\author[B. Bekmaganbetov]{Bekarys Bekmaganbetov}
\address{Division of Applied Mathematics, Brown University, 182 George Street, Providence, RI 02912, USA}
\email{bekarys\_bekmaganbetov@brown.edu}

\author[H. Dong]{Hongjie Dong}
\address{Division of Applied Mathematics, Brown University, 182 George Street, Providence, RI 02912, USA}\email{hongjie\_dong@brown.edu}
\thanks{H. Dong and B. Bekmaganbetov were partially supported by the NSF under agreement DMS-2350129.}


\begin{abstract}
We prove the well-posedness and regularity of solutions in mixed-norm weighted Sobolev spaces for a class of second-order parabolic and elliptic systems in divergence form in the half-space $\R^d_+ = \{x_d > 0\}$ subject to the conormal boundary condition. Our work extends results previously available for scalar equations to the case of systems of equations. The leading coefficients are the product of $x_d^{\alpha}$ and bounded non-degenerate matrices, where $\alpha \in (-1,\infty)$. The leading coefficients are assumed to be merely measurable in the $x_d$ variable, and to have small mean oscillations in small cylinders with respect to the other variables. If the parameter $\alpha>0$, the lower-order coefficients are allowed to blow-up near the boundary. Our results readily generalize to infinite-dimensional equations in general real and complex Hilbert spaces.
\end{abstract}
\maketitle

\section{Introduction}
Let $d,n \in \mathbb{N}, T \in (-\infty, \infty], \R_+ = (0,\infty)$,
\[
\R^d_+ = \{x=(x_1,\dots,x_d) \in \R^d \colon x_d \in \R_+\}
\]
is the upper half-space, $\Omega_T = (-\infty, T) \times \R^d_+$, and $\alpha \in (-1,\infty)$ is a parameter.

We study second-order parabolic systems of $n$ equations in divergence form with singular-degenerate coefficients of the type
 \begin{equation} \label{eq:eqpar}
 \begin{gathered}
 x_d^{\alpha}a_0(x_d) u_t - D_i(x_d^{\alpha}A^{ij}D_j u+x_d^{\alpha-1}B^i u) + x_d^{\alpha-1}\Hat{B}^i D_i u + x_d^{\alpha-2} C u + \lambda x_d^{\alpha} C_0 u  \\
  = -D_i(x_d^{\alpha}F^i) + x_d^{\alpha-1} F_0 + \sqrt{\lambda}x_d^{\alpha} f 
 \end{gathered} \quad \text{in $\Omega_T$}
 \end{equation}
with the conormal boundary condition
\begin{equation} \label{eq:conormal}
\lim_{x_d \to 0^+} x_d^{\alpha}(A^{dj}D_j u+ x_d^{-1}B^d u - F^d) = 0 \quad \text{on $(-\infty,T) \times \R^{d-1}$}.
\end{equation}

Here for $i,j=1,\dots,d$, the coefficients $A^{ij}$, $B^i$, $\Hat{B}^i$, $C$, $C_0:\Omega_T \to \R^{n \times n}$ and $a_0=a_0(x_d):\R_+ \to \R^{n \times n}$ are $n \times n$ matrix-valued measurable functions on $\Omega_T$, $F^i,F_0, f: \Omega_T \to \R^n$ are given (column-) vector-valued measurable functions on $\Omega_T$, and $u:\Omega_T \to \R^n$ is the unknown solution. That is, $u=[u^1,\dots,u^n]^T \in \R^n$, and similarly for $F^i,F_0$, and $f$; $A^{ij}=[A^{ij,k}_l]_{1\leq k,l \leq n} \in \R^{n\times n}$, and similarly for $a_0,B^i,\Hat{B}^i,C$, and $C_0$. We impose standard non-degeneracy, boundedness, and strong ellipticity conditions on $a_0, C_0$, and $A^{ij}$; see (\ref{eq:ellipticity}) and (\ref{eq:nondega0c0}) in Section \ref{sec:2}. We naturally identify $n\times n$ matrices (e.g., the coefficients $a_0,A^{ij},B^i,\Hat{B}^i,C,C_0$) with linear operators $\R^n \to \R^n$. Here $Cu$ stands for the product of a matrix $C$ and a column-vector $u$, that is, $(Cu)^k = \sum_{l=1}^n C^k_l u^l$ for $1 \leq k \leq n$, and similarly for the other terms in (\ref{eq:eqpar}). Moreover, by $D_i = \frac{\partial}{\partial x_i}$ we denote the spatial weak partial derivatives, $\lambda \geq 0$ is a parameter, and we employ the summation convention over repeated upper and lower indices throughout the paper.

Depending on the value of $\alpha$, the coefficients $x_d^{\alpha}a_0$, $x_d^{\alpha}A^{ij}$, $x_d^{\alpha-1}B^i$, $x_d^{\alpha-1}\Hat{B}^i$, and $x_d^{\alpha-2}C$ in (\ref{eq:eqpar}) may be singular or degenerate near the boundary $\partial \R^d_+$. Equations with singular or degenerate coefficients appear naturally in various problems of analysis and applied mathematics. Such equations are closely related to operators of fractional order via the so-called extension
operators, and are therefore an important tool in the study of fractional equations. See \cite{CaffarelliSilvestre2007CPDE, StingaTorrea2017SIAM} for the construction of the extension operators for the fractional Laplace and the fractional heat equations, and \cite{CRS2010JEMS, DeSilvaSavinSire2014} for applications of equations with singular-degenerate coefficients to the study of regularity in free boundary problems associated to the fractional Laplacian. See also \cite{TTV2024ARMA, AFV2025RevMatIber}, where degenerate equations are used to prove the regularity of ratios of solutions to uniformly elliptic equations, i.e., the boundary Harnack principles. We refer the reader to \cite{DongPhan2021TAMS} for more references on connections of equations with singular-degenerate coefficients to various problems in geometric PDEs, porous media, probability, mathematical finance, and mathematical biology.

This paper has two main objectives. First, we consider (\ref{eq:eqpar})-(\ref{eq:conormal}) with lower-order coefficients $B^i=\Hat{B}^i=C=0$, or slightly more generally, satisfying growth bounds
\begin{equation} \label{eq:intro1}
|B| \lesssim x_d, \quad |\Hat{B}| \lesssim x_d, \quad |C| \lesssim x_d^{2}.
\end{equation}
\noindent
We show that under some regularity assumptions on $A^{ij}$ and $C_0$, for any $p,q \in (1,\infty)$ and wide classes of weights $w=w(t,x)=w_0(t)w_1(x)$, holds the estimate
\begin{equation} \label{eq:estintro}
\Vert Du \Vert_{L_{q,p}(\Omega_T, wd\mu)} + \sqrt{\lambda} \Vert u \Vert_{L_{q,p}(\Omega_T, wd\mu)} \leq N(\Vert F \Vert_{L_{q,p}(\Omega_T, wd\mu)} + \Vert f \Vert_{L_{q,p}(\Omega_T, wd\mu)}),
\end{equation}
where $N>0$ does not depend on $u,F$, $f$, and also the dimension $n$. Here $L_{q,p}(\Omega_T, wd\mu)$ is a mixed-norm weighted Lebesgue space with the norm
\[
\Vert f \Vert_{L_{q,p}(\Omega_T, wd\mu)} = \left(\int_{-\infty}^T \left(\int_{\R^d_+}|f(t,x)|^p w_1(x)x_d^{\alpha}dx\right)^{q/p}w_0(t)dt\right)^{1/q},
\]
$d\mu$ stands for $x_d^{\alpha}dxdt$, which is a natural measure for $L_p$ estimates for (\ref{eq:eqpar})-(\ref{eq:conormal}), and we assume that $\lambda>0$ is sufficiently large. In addition to estimate (\ref{eq:estintro}), we also obtain the corresponding existence and uniqueness of solutions. See Theorem \ref{thm:main1} (i) and Section \ref{sec:2} for precise statements and definitions. In the case of scalar equations ($n=1$), the result of our Theorem \ref{thm:main1} (i) was known from the recent work \cite{DongPhanIndiana2023}. The arguments in \cite{DongPhanIndiana2023} cannot be directly generalized to the case of systems due to the use of Moser's iteration, a purely real scalar technique, in one of the key steps of the proof. Showing that estimate (\ref{eq:estintro}) and solvability of (\ref{eq:eqpar})-(\ref{eq:conormal}) still hold in the case of systems with $n \geq 2$ was one of our motivations to write this paper. Moreover, our results hold for solutions taking values in general real or complex Hilbert spaces $H$. That is, if $u,F^i,F_0,f: \Omega_T \to H$ are $H$-valued, and the coefficients $a_0,A^{ij},B^i,\Hat{B}^i,C,C_0:\Omega_T \to \mathcal{B}(H)$, where $\mathcal{B}(H)$ is the space of bounded linear operators on $H$. However, we impose more restricitive measurability conditions on the operator-valued coefficients if the space $H$ is non-separable, see Section \ref{sec:7} for details.

Our second objective in this work is to show that the results described above remain valid under the weaker assumption that the lower-order coefficients $B, \Hat{B}$, and $C$ are bounded and small near the boundary of $\R^d_+$, compared to (\ref{eq:intro1}). That is, we consider
\begin{equation} \label{eq:introlower}
|B| \leq \varepsilon, \quad |\Hat{B}| \leq \varepsilon, \quad |C| \leq \varepsilon
\end{equation}
near the boundary, where $\varepsilon>0$ is sufficiently small. However, the additional assumption $\alpha>0$ is necessary for $|B| \leq \varepsilon$ and $|\Hat{B}| \leq \varepsilon$, and $\alpha>1$ for $|C| \leq \varepsilon$, along with some additional natural conditions on the weights $w$. See Theorem \ref{thm:main1} (ii) and Theorem \ref{thm:mainYp} for precise statements. While conditions of type (\ref{eq:introlower}) are common in the study of equations with the {\it zero Dirichlet boundary condition}, see \cite{DongKim2015AdvMath, KimKimLee2022CPAA, KimKHKimLeeJMAA2014}, to the best of our knowledge, our paper is the first one in which the coefficients $B,\Hat{B}$, and $C$ are allowed to satisfy (\ref{eq:introlower}) with the {\it conormal boundary condition}.

Let us point out that the degeneracy of the coefficients provided by $\alpha>0$ and $\alpha>1$ is crucial in order to allow (\ref{eq:introlower}) in the conormal boundary condition setting. One can directly see that some degeneracy of the coefficients is needed just to define weak solutions to (\ref{eq:eqpar})-(\ref{eq:conormal}); see (\ref{eq:soldefvar}) in Section \ref{sec:2} for the definition of weak solutions. Indeed, let $u$ be a solution to (\ref{eq:eqpar})-(\ref{eq:conormal}) and let $|B|,|\Hat{B}|$, and $|C|$ be comparable to a small constant near the boundary $\partial \R^d_+$. Then $x_d^{\alpha-1}|u|$, $x_d^{\alpha-1}|Du|$, and $x_d^{\alpha-2}|u|$, respectively, need to be integrable near the boundary $\partial \R^d_+$. One then immediately notices that even if $u,Du \in L^{\infty}$ up to the boundary, conditions $\alpha>0$ and $\alpha>1$ become necessary.

We now briefly discuss the regularity assumptions on the coefficients imposed in this paper, apart from the standard ones (\ref{eq:ellipticity})-(\ref{eq:nondega0c0}). We assume that all $a_0$, $A^{ij}$, and $C_0$ are merely measurable in the $x_d$ variable and $a_0=a_0(x_d)$ depends on $x_d$ only. Equations involving coefficients without any regularity in one variable arise in problems of linearly elastic laminates and composite materials, for example, in homogenization of layered materials (see, e.g., \cite{ChKV-C1986ARMA}). Furthermore, we assume that $A^{ij}$ and $C_0$ are partially VMO, that is, they have small mean oscillations on small cylinders in the $t,x_1,\dots,x_{d-1}$ variables. See Assumption \ref{assump} below. The class of partially VMO coefficients was first introduced in \cite{KimKrylov2007SIAM, KimKrylov2007PotAn} and has since been widely used in the regularity theory for elliptic and parabolic PDEs.

We now review some more results from the literature related to (\ref{eq:eqpar}), (\ref{eq:conormal}), and results addressing unbounded lower-order coefficients. In \cite{DongPhan2021CalcVar}, the authors studied systems (\ref{eq:eqpar})-(\ref{eq:conormal}) with more general factors $\mu(x_d)$ in place of $x_d^{\alpha}$, and obtained (\ref{eq:estintro}) in the unmixed and unweighted setting, that is, for $q=p \in (1,\infty)$ and $w=1$. There, $\mu(x_d)$ are Muckenhoupt weights of class $A_2(\R)$ satisfying a certain local $A_1$ type condition that cover $\mu(x_d)=x_d^{\alpha}$ in the range $\alpha \in (-1,1)$. We recall that Harnack inequalities and H\"older regularity for equations with general $A_2$ weights as factors in coefficients were obtained in \cite{FabesKenigSerapioni1982CPDE, ChiarenzaSerapioni1985Padova}. Note that one may consider systems analogous to (\ref{eq:eqpar})-(\ref{eq:conormal}) on domains $\Omega$ with curved boundaries by changing $x_d^{\alpha}$ to $\operatorname{dist}(x,\partial \Omega)^{\alpha}$. We will address such systems on bounded domains in a subsequent paper. We also refer to recent work \cite{DongPhanSire2024JGeomAn}, where the authors obtained $L_p$ estimates for general quasilinear elliptic equations on bounded Lipchitz domains with small Lipschitz constants.

In \cite{DongPhan2021RevMatIber} and \cite{DongPhanIndiana2023} the authors considered scalar ($n=1$) non-divergence form analogues of (\ref{eq:eqpar})-(\ref{eq:conormal}) in the ranges $\alpha \in (-1,1)$ and $\alpha \in (-1,\infty)$, respectively. There, an additional structural condition $A^{dj}/A^{dd} \equiv \text{const}$ for $1 \leq j \leq d-1$ is imposed. The results in \cite{DongPhan2021RevMatIber, DongPhanIndiana2023} essentially use the scalar structure of the equations and are not readily generalized to systems of equations. Moreover, it remains an open problem whether classical $W^2_p$ estimates hold for uniformly elliptic and parabolic systems with partially VMO coefficients. We therefore do not consider systems in non-divergence form in this paper and leave this question for future investigation.

Recent results on boundary Schauder-type estimates for parabolic equations with the conormal boundary condition can be found in \cite{AFV2024CVPDE, AFV2025RevMatIber, DongJeon2025arxiv}, and for elliptic equations in \cite{DJV2024CVPDE, TTV2024ARMA}. In these works, Schauder estimates for equations with singular-degenerate coefficients are used to derive boundary Harnack principles for equations with bounded uniformly elliptic coefficients. Parabolic and elliptic systems in divergence form with the zero Dirichlet boundary condition were studied in \cite{DongPhan2021TAMS}, whereas scalar equations in non-divergence form were considered in \cite{DongPhan2023JFA}. In the recent work \cite{BekmaganbetovDong2024arxiv}, the authors studied the solvability of equations in divergence form with the inhomogeneous Dirichlet boundary data taken in spaces of distributions with negative order of differentiability. In these papers, the natural range of the parameter $\alpha$ for equations with the Dirichlet boundary condition is $\alpha \in (-\infty, 1)$.

Various phenomena in pure and applied mathematics are modeled by equations whose coefficients in front of $u_t, D^2u, Du$, and $u$ may exhibit different rates of decay or blow-up. Examples of such anisotropy in coefficients include the Black-Scholes-Merton equation in mathematical finance, degenerate viscous Hamilton-Jacobi equations, and equations involving Baouendi-Grushin operators arising in sub-Riemannian geometry. For recent results on weighted $L_p$ estimates for such anisotropic equations we refer the reader to \cite{DongRyu2026TAMS, DongRyu2025arxiv, DongPhanTranTAMS2023, DongPhanTranJFA2024, MNS2023JDE, FangPhan2510.20051}.

$L_p$ estimates for elliptic and parabolic equations with unbounded lower-order terms is an extensively studied topic, particularly in the classical case $\alpha=0$. Such results are often obtained under some integrability assumptions on the coefficients $B, \Hat{B}$, and $C$; common conditions include $B, \Hat{B} \in L_d$, $C \in L_{d/2}$, or their generalizations. Under these assumptions, $L_p$ estimates can be established for both Dirichlet and conormal boundary conditions. For detailed results in this direction, see \cite{KangKim2017CPAA, Krylov2021TAMS, KimKwon2022TAMS} for elliptic equations and \cite{KimRyuWoo2022JEvolEq} for parabolic equations. As we mentioned earlier, $L_p$ estimates for equations with the zero Dirichlet boundary condition under assumptions of type (\ref{eq:introlower}) are also well studied in the literature; see \cite{DongKim2015AdvMath, KimKimLee2022CPAA, KimKHKimLeeJMAA2014}.

Lastly, we note that the literature on infinite-dimensional parabolic equations in general Banach spaces, including Hilbert spaces, is highly active and extensive. We refer the reader to \cite{FT2002AnnProb} for semilinear stochastic evolution equations in Hilbert spaces with applications to stochastic optimal control and Hamilton-Jacobi-Bellman equations, \cite{Masiero2008SIAMJOptCtrl} for stochastic optimal control problems in Banach spaces; see also \cite{vNVW2008JFA, vNVW2012AnnProb} and references therein. For deterministic equations, numerous maximal $L^p-L^q$ regularity results are available for equations in the so called UMD Banach spaces, an example of which are Hilbert spaces. See \cite{DHP2007MathZ} for general parabolic initial-boundary value problems, \cite{Lindemulder2020JEvolEq, LV2020JDE} for generalizations to weighted spaces, and references therein. We note that parabolic equations in UMD spaces are usually studied via an approach quite different from the one in the present paper. There, parabolic equations are usually viewed as abstract evolution equations and studied via functional analytic tools such as the theory of sectorial operators and $H^{\infty}$-calculus.

In our proofs, we first consider systems without lower-order terms, i.e., when $B=\Hat{B}=C=0$. In this case, our arguments rely on mean oscillation estimates together with Fefferman-Stein type inequalities. We begin with systems having simple coefficients, namely when $A^{ij}$ and $C_0$ depend on $x_d$ only. The key intermediate result here is Lemma \ref{lem:bdrbddnssest}, where we establish several boundary $L_{\infty}$ estimates for solutions of the corresponding homogeneous systems. For scalar equations ($n=1$), such estimates were previously obtained in \cite{DongPhanIndiana2023} via Moser's iteration combined with a certain bootstrap argument. Here we develop a different bootstrap scheme that is applicable to systems of equations and avoids scalar-specific techniques. By combining these $L_{\infty}$ estimates with $L_2$ energy estimates, we derive Lipschitz and mean oscillation bounds for $D_{x'}u$ and $\mathcal{U}=A^{dj}D_j u$. We then prove the solvability of systems with simple coefficients in unweighted $L_p$ spaces without mixed norms. For the general case when $A^{ij}$ and $C_0$ may depend on all the $t,x$ variables, we apply a perturbation argument together with the generalized version of the Fefferman-Stein inequality for weighted mixed-norm spaces established in \cite{DongKimTAMS2018}. Lower-order terms involving $B,\Hat{B},$ and $C$ under condition (\ref{eq:intro1}) are then added by taking $\lambda$ sufficiently large. In the singular case when the lower-order terms satisfy (\ref{eq:introlower}), they are treated as perturbations and controlled using certain weighted Hardy-type inequalities. We use the results from \cite{Muckenhoupt1972} to identify the additional conditions on weights $w$ under which such Hardy-type inequalities hold.

In our proofs, we do not rely on any specific properties of $\R^n$ other than its Hilbert space structure, nor do we use the base field being $\R$ in any essential way. Since our initial motivation was studying systems of equations, we chose to write the proofs in the language of systems. We briefly explain in Section \ref{sec:7} how our results generalize to equations in general real or complex Hilbert spaces. The Hilbert space structure, however, is essential for deriving the $L_2$ estimates in Section \ref{sec:L2}. It is of interest to understand whether our results can be generalized to equations in more general spaces, for example, UMD spaces.

The rest of the paper is organized as follows. In Section \ref{sec:2} we introduce the notation, function spaces, and state the main results, Theorems \ref{thm:main1} and \ref{thm:mainYp}. In Section \ref{sec:L2} we prove the existence and uniqueness of $L_2$ solutions. In Section \ref{sec:4} we consider systems having simple coefficients. In Section \ref{sec:5} we prove Part (i) of Theorem \ref{thm:main1}, where $B, \Hat{B}$, and $C$ satisfy (\ref{eq:intro1}). Section \ref{sec:lowerterms} addresses systems with singular lower-order terms satisfying (\ref{eq:introlower}), completing the proofs of Theorems \ref{thm:main1} and \ref{thm:mainYp}. Finally, in Section \ref{sec:7} we discuss the generalization of our results to Hilbert space-valued solutions.

\section{Notation and main results} \label{sec:2}

\subsection{Notation} We equip $\R^n$ with the standard dot product $a\cdot b := \sum_{k=1}^n a^k b^k$.

For $p \in (1,\infty)$, we denote by $p'=\frac{p}{p-1}$ the H\"older dual to $p$.
By $N_{a,b,\dots}$ we denote positive constants which depend on the subscripted indices $a,b,\dots$.

For $r > 0$ and $z_0=(t_0,x_0)$ with $x_0 = (x_0', x_{0d}) \in \R^d$ and $t_0 \in \R$, we define $z_0' = (t_0, x_0')$, 
\begin{equation*}
\begin{gathered}
B_r'(x_0') = \{x' \in \R^{d-1}: |x'-x_0'| < r\}, \quad B_r(x_0) = \{x \in \R^d: |x-x_0| < r\}, \\
D_r(x_0) = B_r'(x_0') \times (x_{0d}-r, x_{0d} + r), \quad D_r^+(x_0) = D_r(x_0) \cap \R^d_+, \\
Q_r(z_0) = (t_0-r^2, t_0) \times D_r(x_0), \quad Q_r^+(z_0) = (t_0-r^2, t_0) \times D_r^+(x_0), \\
Q_r'(z_0') = (t_0-r^2, t_0) \times B_r'(x_0').
\end{gathered} 
\end{equation*}
When $t_0 = 0$ and $x_0 = 0$, we denote $B_r' = B_r'(0)$, $B_r = B_r(0)$, $D_r = D_r(0)$, $D_r^+ = D_r^+(0)$, $Q_r = Q_r(0)$, $Q_r^+ = Q_r^+(0)$, and $Q_r' = Q_r'(0)$.

By $D_i u = \frac{\partial u}{\partial x_i}$, $i=1,\dots,d$, we denote the spatial weak partial derivatives of $u$, $Du = (D_1 u, \dots, D_d u)$, and $D_{x'}u = (D_1 u, \dots, D_{d-1} u)$.

Throughout the paper, for $A=[A^{ij}]^{1\leq i,j \leq d}$, $A^{ij} \in \R^{n\times n}$, we denote
\[
|A| = \sup_{\xi,\eta \in (\R^{n})^d, |\xi|=|\eta|=1} \left|\sum_{i,j=1}^d A^{ij}\xi_j \cdot \eta_i\right|,
\]
where for $\xi = (\xi_1,\dots,\xi_d) \in (\R^{n})^d$ we set $|\xi|^2 = \sum\limits_{i=1}^d |\xi_i|^2$, and similarly for $\eta$. In particular, $|A^{ij}\xi_j \cdot \eta_i| \leq |A||\xi||\eta|$ for all $\xi,\eta \in (\R^{n})^d$. For a matrix $C \in \R^{n \times n}$, by $|C|$ we denote its operator norm
\[
|C| := |C|_{\text{op}} = \sup_{u \in \R^n, |u|=1} |Cu|.
\]
For $B=(B^1,\dots,B^d)$ and $\Hat{B}=(\Hat{B}^1,\dots,\Hat{B}^d)$, where $B^i, \Hat{B}^i \in \R^{n\times n}$, we set
\[
|B|_c^2 = \left|\sum_{i=1}^d (B^i)^T B^i\right|_{\text{op}}, \quad |\Hat{B}|_r^2 = \left|\sum_{i=1}^d \Hat{B}^i (\Hat{B}^i)^T\right|_{\text{op}}.
\]
Here $|\cdot|_c$ stands for ``column" operator norm and $|\cdot|_r$ stands for ``row" operator norm. In particular, for all $u,\xi_1,\dots,\xi_d \in \R^n$ it holds that
\[
|(B^1u,\dots,B^d u)| = \sqrt{\sum_{i=1}^d |B^i u|^2} \leq |B|_c|u|, \quad \left|\sum_{i=1}^d \Hat{B}^i \xi_i \right| \leq |\Hat{B}|_r|\xi|.
\]
Moreover, $|B|_c$ and $|\Hat{B}|_r$ are the best constants in the inequalities above. To simplify notation, we will omit the subscripts in $|B|_c$ and $|\Hat{B}|_r$ and simply write $|B|$ and $|\Hat{B}|$ throughout the paper. We understand $|\Hat{B}|$ as $|\Hat{B}|_r$ if the ``hat" symbol is present, and $|B|=|B|_c$ otherwise.

We assume that the leading coefficients $A^{ij}$ are bounded and satisfy the strong ellipticity condition as follows: there exists $\kappa \in (0,1]$ such that
\begin{equation} \label{eq:ellipticity}
\kappa \sum_{i=1}^d |\xi_i|^2 \leq \sum_{i,j=1}^d A^{ij}\xi_j \cdot \xi_i, \quad |A| \leq \kappa^{-1},
\end{equation}
on $\Omega_T$ for all $\xi_i \in \R^{n}, i = 1, \dots, d$.  The coefficient matrices $a_0$ and $C_0$ are assumed to be bounded and positive definite as follows: 
\begin{equation} \label{eq:nondega0c0}
\begin{gathered}
 \kappa |\xi|^2 \leq a_0 \xi \cdot \xi, \quad |a_0| \leq \kappa^{-1}, \quad a_0 = a_0^T, \\
 \kappa |\xi|^2 \leq C_0 \xi \cdot \xi, \quad |C_0| \leq \kappa^{-1},
\end{gathered}
\end{equation}
on $\Omega_T$ for all $\xi \in \R^n$.

Let $\mathcal{O} \subset \R^d_+$ denote an open set of the form $\mathcal{O} = \mathcal{O}' \times (a,b)$, where $\mathcal{O}' \subset \R^{d-1}$ is open and $0 \leq a < b \leq \infty$. For such $\mathcal{O}$, we set $\mathcal{O}_{\text{test}} = \mathcal{O}' \times (-\infty,b)$ if $a = 0$, and $\mathcal{O}_{\text{test}} = \mathcal{O}$ if $a>0$. Also let $-\infty \leq S < T \leq +\infty$.

We denote $d\mu = x_d^{\alpha} dxdt$, $d\mu=x_d^{\alpha} dx$, or $d\mu = x_d^{\alpha} dx_d$. It will always be clear from context which measure is used in a particular situation. We say that $u$ is a weak solution to equation (\ref{eq:eqpar}) in $(S,T) \times \mathcal{O}$ if
\begin{equation} \label{eq:soldefvar}
\begin{gathered}
\int_{(S,T) \times \mathcal{O}}(-a_0 u \cdot \varphi_t + (A^{ij}D_j u + x_d^{-1} B^i u) \cdot D_i \varphi + x_d^{-1} \Hat{B}^i D_i u \cdot \varphi + x_d^{-2} C u \cdot \varphi + \lambda C_0 u \cdot \varphi) d\mu \\
= \int_{(S,T) \times \mathcal{O}}(F^i \cdot D_i \varphi + x_d^{-1} F_0 \cdot \varphi + \sqrt{\lambda} f\cdot \varphi) d\mu
\end{gathered}
\end{equation}
holds for all $\varphi \in C_0^{\infty}((S,T) \times \mathcal{O}; \R^n)$. If $a=0$, we say that $u$ is a weak solution to equation (\ref{eq:eqpar}) in $(S,T) \times \mathcal{O}$ with the conormal boundary condition (\ref{eq:conormal}) satisfied on $(S,T) \times \mathcal{O}'$ if (\ref{eq:soldefvar}) holds for all $\varphi \in C_0^{\infty}((S,T) \times \mathcal{O}_{\text{test}}; \R^n)$. For the integrals in (\ref{eq:soldefvar}) to be well defined, we will always have
\begin{equation*}
x_d^{\alpha} (|u|+|Du|+x_d^{-1}|B| |u|+x_d^{-1}|\Hat{B}| |D u|+x_d^{-2} |C| |u|+|F|+x_d^{-1}|F_0|+|f|) \in L_{1,\text{loc}}\left((S,T)\times \left(\mathcal{O}_{\text{test}} \cap \overline{\R^d_+}\right)\right).
\end{equation*}

For $p, q \in (1, \infty)$, weights (that is, measurable and almost everywhere positive finite functions) $w_0(t)$ and $w_1(x)$, and a Banach space $X$, we define $L_{q,p}((S,T) \times \mathcal{O}, wd\mu; X)$ to be the weighted mixed-norm Lebesgue space of measurable (strongly measurable in the sense of Bochner if $\operatorname{dim}X=\infty$) $X$-valued functions with the norm
\[
\Vert f \Vert_{L_{q,p}((S,T) \times \mathcal{O}, wd\mu; X)} = \left(\int_S^T \left( \int_{\mathcal{O}} |f(t,x)|_X^p w_1(x) d\mu(x) \right)^{q/p} w_0(t)dt \right)^{1/q}.
\]
If $X=\R$ or $X=\C$, we write $L_{q,p}(\cdot)$ instead of $L_{q,p}(\cdot;X)$. When the space $X$ is understood from context, we will often write $f \in L_{q,p}(\cdot)$ meaning $|f|_X \in L_{q,p}(\cdot)$. 

In general, if $(v_{k})$ is a finite collection of vectors $v_k \in \R^n$, we set $|(v_k)|^2 = \sum_k |v_k|^2$. For instance, $|Du|^2=\sum_{i=1}^d |D_i u|^2$, $|(F^1,\dots,F^d)|^2 = \sum_{i=1}^d |F^i|^2$, etc.

We define
\begin{gather}
\mathbb{H}^{-1}_{q,p}((S,T)\times \mathcal{O}, wd\mu; \R^n) = \big\{g \in \mathcal{D}'((S,T)\times\mathcal{O}_{\text{test}}; \R^n): \; \langle g, \varphi \rangle = \int_{(S,T)\times\mathcal{O}}(F^i \cdot D_i \varphi + f \cdot \varphi)d\mu \nonumber \\
\text{for all} \; \varphi \in C_0^{\infty}((S,T)\times\mathcal{O}_{\text{test}}; \R^n) \; \text{for some $F^1,\dots,F^d,f \in L_{q,p}((S,T)\times \mathcal{O}, wd\mu; \R^n)$}\big\}, \label{eq:H-1vardef} \\
\Vert g \Vert_{\mathbb{H}^{-1}_{q,p}((S,T)\times \mathcal{O}, wd\mu; \R^n)} = \inf \left\{ \Vert F \Vert_{L_{q,p}((S,T)\times \mathcal{O}, wd\mu)} + \Vert f \Vert_{L_{q,p}((S,T)\times \mathcal{O}, wd\mu)} : \eqref{eq:H-1vardef} \; \text{holds}\right\}, \nonumber\\
\mathcal{H}^{1}_{q,p}((S,T)\times \mathcal{O}, wd\mu; \R^n) = \bigl\{u \in L_{q,p}((S,T)\times \mathcal{O}, wd\mu; \R^n): \nonumber \\
Du \in L_{q,p}((S,T)\times \mathcal{O}, wd\mu), x_d^{\alpha}a_0(x_d) u_t \in \mathbb{H}^{-1}_{q,p}((S,T)\times \mathcal{O}, wd\mu; \R^n)\bigr\}, \nonumber \\
\Vert u \Vert_{\mathcal{H}^{1}_{q,p}((S,T)\times \mathcal{O}, wd\mu; \R^n)}  \nonumber \\ = \Vert Du \Vert_{L_{q,p}((S,T)\times \mathcal{O}, wd\mu)} + \Vert u \Vert_{L_{q,p}((S,T)\times \mathcal{O}, wd\mu)} + \Vert x_d^{\alpha}a_0(x_d)u_t \Vert_{\mathbb{H}^{-1}_{q,p}((S,T)\times \mathcal{O}, wd\mu; \R^n)}. \nonumber
\end{gather}
Here we understand $x_d^{\alpha} a_0(x_d)u_t$ as an element of $\mathcal{D}'((S,T)\times \mathcal{O}_{\text{test}}; \R^n)$ given by
\[
\langle x_d^{\alpha} a_0(x_d)u_t, \varphi \rangle = -\int_{(S,T)\times \mathcal{O}}a_0(x_d)u \cdot \varphi_t d\mu
\]
for $\varphi \in C_0^{\infty}((S,T)\times \mathcal{O}_{\text{test}}; \R^n)$. Note that the definition of the space $\mathcal{H}^1_{q,p}$ depends on the choice of the coefficient matrix $a_0(x_d)$.

When $p=q$, we write $L_p((S,T) \times \mathcal{O}, wd\mu)$ instead of $L_{p,p}((S,T) \times \mathcal{O}, wd\mu)$ and similarly for the other function spaces. When both $w_0 \equiv 1$ and $w_1 \equiv 1$, we have $L_{q,p}((S,T) \times \mathcal{O}, wd\mu) = L_{q,p}((S,T) \times \mathcal{O}, \mu)$, and similarly for the other function spaces.

In the time-independent case, for $w=w(x)$ we set
\begin{gather*}
\Vert f \Vert_{L_p(\mathcal{O}, w d\mu)} = \left(\int_{\mathcal{O}} |f(x)|^p w(x)d\mu(x)\right)^{1/p}, \\
\Vert u \Vert_{W^1_p(\mathcal{O}, w d\mu)} = \Vert Du \Vert_{L_p(\mathcal{O}, w d\mu)} + \Vert u \Vert_{L_p(\mathcal{O}, w d\mu)}.
\end{gather*}

For a set $\Omega$, a measure $\nu$ on $\Omega$ with $0 < \nu(\Omega) < \infty$, and a function $f \in L_1(\Omega, \nu)$, we write
\[
(f)_{\Omega,\nu} = \fint_{\Omega} f d\nu = \frac{1}{\nu(\Omega)} \int_{\Omega} f d\nu.
\]

For $p \in (1,\infty)$, a weight $w$ on $\R^d_+$ is said to be in the $A_p(\R^d_+, \mu)$ Muckenhoupt class if
\[
[w]_{A_p(\R^d_+,\mu)} = \sup_{r>0,x\in \overline{\R^d_+}} \left(\fint_{D_r^+(x)}w d\mu\right) \times \left(\fint_{D_r^+(x)}w^{-\frac{1}{p-1}} d\mu\right)^{p-1} < \infty.
\]

Similarly, a weight $w$ on $\R^m$ belongs to the class $A_p(\R^m)$ if
\[
[w]_{A_p(\R^m)} = \sup_{r>0, y_0 \in \R^m} \left(\fint_{B_r(y_0)} w(y) dy\right) \times \left(\fint_{B_r(y_0)} w^{-\frac{1}{p-1}}(y) dy\right)^{p-1} < \infty.
\]

We will also need the following two classes of weights $X_p(\mu)$ and $Y_p(\mu)$.
\begin{defn} 
Let $\alpha \in \R, d\mu(y) = y^{\alpha}dy$ for $y \in \R_+$, and $p \in (1,\infty)$.

We say that a weight $w$ on $\R_+$ belongs to the class $X_p(\mu)$ if
\[
[w]_{X_p(\mu)} = \sup_{r>0} \left[\int_0^r w(y) y^{-p} d\mu(y)\right]^{\frac{1}{p}} \left[ \int_r^{\infty}w(y)^{-\frac{1}{p-1}}y^{-\alpha p'} d\mu(y)\right]^{\frac{1}{p'}} < \infty.
\]

We say that a weight $w$ on $\R_+$ belongs to the class $Y_p(\mu)$ if
\[
[w]_{Y_p(\mu)} = \sup_{r>0} \left[\int_r^{\infty} w(y) y^{-\alpha p} d\mu(y)\right]^{\frac{1}{p}} \left[ \int_0^r w(y)^{-\frac{1}{p-1}}y^{-p'} d\mu(y)\right]^{\frac{1}{p'}} < \infty.
\]
\end{defn}

\subsection{Main theorems}

We impose the following partially BMO regularity assumption on the coefficents $A^{ij}$ and $C_0$. Let $\gamma_0 > 0$ and $R_0 \in (0,\infty]$.

\begin{assumption}[$\gamma_0, R_0$] \label{assump}
For any $r \in (0,R_0]$ ($r \in (0,\infty)$ if $R_0 = \infty$) and $z_0=(z_0',x_{0d}) \in \overline{\Omega_T}$, there exist measurable $[A^{ij}]_{r,z_0}(x_d),[C_{0}]_{r,z_0}(x_d):\R_+\to \R^{n\times n}$, $1 \leq i,j \leq d$, satisfying (\ref{eq:ellipticity})-(\ref{eq:nondega0c0}) for a.e. $x_d \in \R_+$, and
\[
\fint_{Q_r^+(z_0)}|A(t,x)-[A]_{r,z_0}(x_d)| d\mu(z) \leq \gamma_0, \quad \fint_{Q_r^+(z_0)}|C_{0}(t,x)-[C_{0}]_{r,z_0}(x_d)| d\mu(z) \leq \gamma_0,
\]
where $d\mu(z) = x_d^{\alpha} dxdt$.
\end{assumption}

\begin{thm} \label{thm:main1}
Let $\alpha \in (-1, \infty), p,q \in (1,\infty), \kappa \in (0,1], K \in [1,\infty), R_0 \in (0,\infty]$, and $K_1 \in [0,\infty)$. Let $P=(d,\alpha,p,q,\kappa,K)$ denote the collection of listed parameters. There exist $\gamma_0=\gamma_0(P)>0$ and $\lambda_0=\lambda_0(P)>0$ such that the following assertions hold. Let $w = w_0(t)w_1(x)$, where $w_0 \in A_q(\R), w_1 \in A_p(\R^d_+, \mu)$ with $[w_0]_{A_q(\R)} \leq K, [w_1]_{A_p(\R^d_+,\mu)} \leq K$. Suppose that \eqref{eq:ellipticity}, \eqref{eq:nondega0c0}, and Assumption {\rm \ref{assump}} $(\gamma_0,R_0)$ are satisfied. Also assume that $F_0=0$.

\begin{enumerate}
\item[(i)]
Suppose that $|B| \leq K_1x_d$, $|\Hat{B}| \leq K_1x_d$, and $|C| \leq K_1^2 x_d^2$. There exist $M=M(P)>0$ and $N=N(P)>0$ such that for any $F=(F^1,\dots,F^d) \in L_{q,p}(\Omega_T,w d\mu; \R^n)^{d}$, $f \in L_{q,p}(\Omega_T,w d\mu; \R^n)$, $u \in \mathcal{H}^1_{q,p}(\Omega_T,w d\mu; \R^n)$, and $\lambda \geq \max(\lambda_0 R_0^{-2}, MK_1^2)$ satisfying \eqref{eq:eqpar}-\eqref{eq:conormal}, it holds that
\begin{equation} \label{eq:estmain1}
\Vert Du \Vert_{L_{q,p}(\Omega_T,w d\mu)} + \sqrt{\lambda}\Vert u \Vert_{L_{q,p}(\Omega_T,w d\mu)} \leq N \Vert F\Vert_{L_{q,p}(\Omega_T,w d\mu)} + N \Vert f\Vert_{L_{q,p}(\Omega_T,w d\mu)},
\end{equation}
Moreover, for any $F\in L_{q,p}(\Omega_T,w d\mu; \R^n)^{d}$, $f \in L_{q,p}(\Omega_T,w d\mu; \R^n)$, and $\lambda \geq \max(\lambda_0 R_0^{-2}, MK_1^2)$ with $\lambda > 0$, there exists a unique solution $u \in \mathcal{H}^1_{q,p}(\Omega_T,w d\mu; \R^n)$ to \eqref{eq:eqpar}-\eqref{eq:conormal}.
\item[(ii)]
Assume, additionally, that $\alpha \in (0,\infty)$, and $w_1(x) = w_2(x')w_3(x_d)$, where $w_2 \in A_p(\R^{d-1}), w_3 \in A_p(\R_+, \mu) \cap X_p(\mu)$ with $[w_2]_{A_p(\R^{d-1})}\cdot [w_3]_{A_p(\R_+,\mu)} \leq K, [w_3]_{X_p(\mu)} \leq K_2, K_2 \in (0,\infty)$. There exist $\varepsilon=\varepsilon(P) > 0, M=M(P)>0$, and $N=N(P)>0$, such that the result of {\rm (i)} holds under the assumption
$$
|B| \leq \max(\varepsilon K_2^{-1}, K_1 x_d),\quad 
|\Hat{B}| \leq K_1 x_d,\quad \text{and}\quad
|C| \leq \max(K_1 K_2^{-1} x_d, K_1^2x_d^2).
$$ 
Moreover, in place of \eqref{eq:estmain1} holds the estimate
\begin{equation}\label{eq:estmainXp}
\begin{aligned} 
\Vert Du \Vert_{L_{q,p}(\Omega_T,w d\mu)} + K_2^{-1}\Vert u/x_d \Vert_{L_{q,p}(\Omega_T,w d\mu)} &+ \sqrt{\lambda}\Vert u \Vert_{L_{q,p}(\Omega_T,w d\mu)} \\
 &\leq N \Vert F\Vert_{L_{q,p}(\Omega_T,w d\mu)} + N \Vert f\Vert_{L_{q,p}(\Omega_T,w d\mu)}.
\end{aligned}
\end{equation}
\end{enumerate}
\end{thm}

\begin{thm} \label{thm:mainYp}
Let $\alpha \in (0,\infty), p,q \in (1,\infty), \kappa \in (0,1], K \in [1,\infty), R_0 \in (0,\infty], K_1 \in [0,\infty)$. Let $P=(d,\alpha,p,q,\kappa,K)$ denote the collection of listed parameters. Let $\gamma_0=\gamma_0(P)>0$ and $\lambda_0=\lambda_0(P)>0$ be from Theorem {\rm \ref{thm:main1}}. Let $w = w_0(t)w_2(x')w_3(x_d)$, where $w_0 \in A_q(\R), w_2 \in A_p(\R^{d-1}), w_3 \in A_p(\R_+, \mu) \cap Y_p(\mu)$ with $[w_0]_{A_q(\R)} \leq K, [w_2]_{A_p(\R^{d-1})}\cdot [w_3]_{A_p(\R_+,\mu)} \leq K, [w_3]_{Y_p(\mu)} \leq K_3, K_3 \in (0,\infty)$. Suppose that \eqref{eq:ellipticity}, \eqref{eq:nondega0c0}, and Assumption {\rm \ref{assump}} $(\gamma_0,R_0)$ are satisfied.

\begin{enumerate}
\item[(i)]
There exist $\varepsilon=\varepsilon(P) > 0, M=M(P)>0$, and $N=N(P)>0$ such that the following assertions holds. Suppose that $|B| \leq K_1x_d$, $|\Hat{B}| \leq \max(\varepsilon K_3^{-1}, K_1x_d)$, and $|C| \leq \max(K_1 K_3^{-1} x_d, K_1^2 x_d^2)$. Then for any $F=(F^1,\dots,F^d) \in L_{q,p}(\Omega_T,w d\mu; \R^n)^{d}$, $F_0,f \in L_{q,p}(\Omega_T,w d\mu; \R^n)$, $u \in \mathcal{H}^1_{q,p}(\Omega_T,w d\mu; \R^n)$, and $\lambda \geq \max(\lambda_0R_0^{-2}, MK_1^2)$ satisfying \eqref{eq:eqpar}-\eqref{eq:conormal}, it holds that
\begin{equation} \label{eq:estmainYp}
\begin{aligned}
&\Vert Du \Vert_{L_{q,p}(\Omega_T,w d\mu)} + \sqrt{\lambda}\Vert u \Vert_{L_{q,p}(\Omega_T,w d\mu)} \\
&\leq N \Vert F\Vert_{L_{q,p}(\Omega_T,w d\mu)} + N \Vert f\Vert_{L_{q,p}(\Omega_T,w d\mu)} + N K_3 \Vert F_0 \Vert_{L_{q,p}(\Omega_T,w d\mu)}.
\end{aligned}
\end{equation}
Moreover, for any $F\in L_{q,p}(\Omega_T,w d\mu; \R^n)^{d}$, $F_0,f \in L_{q,p}(\Omega_T,w d\mu; \R^n)$, and $\lambda \geq \max(\lambda_0R_0^{-2}, MK_1^2)$ with $\lambda > 0$, there exists a unique solution $u \in \mathcal{H}^1_{q,p}(\Omega_T,w d\mu; \R^n)$ to \eqref{eq:eqpar}-\eqref{eq:conormal}.
\item[(ii)]
Assume, additionally, that $\alpha \in (1,\infty)$ and $w_3 \in A_p(\R_+, \mu) \cap X_p(\mu) \cap Y_p(\mu)$ with $[w_3]_{X_p(\mu)} \leq K_2, K_2 \in (0,\infty)$. There exist $\varepsilon=\varepsilon(P) > 0, M=M(P)>0$, and $N=N(P)>0$ , such that the result of {\rm (i)} holds under the assumption
\[
|B| \leq \max(\varepsilon K_2^{-1}, K_1x_d), |\Hat{B}| \leq \max(\varepsilon K_3^{-1}, K_1x_d), |C| \leq \max(\varepsilon (K_2 K_3)^{-1}, K_1 K_2^{-1}x_d, K_1 K_3^{-1} x_d, K_1^2 x_d^2).
\]
Moreover, in place of \eqref{eq:estmainYp} holds the estimate
\begin{equation}\label{eq:estmainXpYp}
\begin{aligned} 
\Vert Du \Vert_{L_{q,p}(\Omega_T,w d\mu)} &+ K_2^{-1}\Vert u/x_d \Vert_{L_{q,p}(\Omega_T,w d\mu)} + \sqrt{\lambda}\Vert u \Vert_{L_{q,p}(\Omega_T,w d\mu)} \\
 &\leq N \Vert F\Vert_{L_{q,p}(\Omega_T,w d\mu)} + N \Vert f\Vert_{L_{q,p}(\Omega_T,w d\mu)}+ N K_3 \Vert F_0 \Vert_{L_{q,p}(\Omega_T,w d\mu)}.
\end{aligned}
\end{equation}
\end{enumerate}

\end{thm}

Here in estimates in Theorems \ref{thm:main1} and \ref{thm:mainYp} involving parameters $R_0,K_1,K_2,K_3$, we write the explicit dependence on them in order for the constants $\lambda_0,\varepsilon, M$, and $N$ to be independent of these parameters. Moreover, in the above form, all estimates in Theorems \ref{thm:main1} and \ref{thm:mainYp} are invariant with respect to the parabolic scaling $(t,x) \to (s^2t,sx)$, where $s>0$.

\begin{remark} \label{rem:1}
It is worth highlighting the following few observations about the $X_p$ and $Y_p$ conditions and their roles in Theorems \ref{thm:main1} and \ref{thm:mainYp}.
\begin{enumerate}
\item As a typical example, in Theorems \ref{thm:main1} and \ref{thm:mainYp} we can consider power weights $w_3(y) = y^{\gamma}$, where $\gamma \in \R$. A straightforward calculation shows that
\begin{gather*}
y^{\gamma} \in A_p(\R_+, \mu) \iff \frac{1+\alpha+\gamma}{p} \in (0, 1+\alpha),  \\
y^{\gamma} \in X_p(\mu) \iff \frac{1+\alpha+\gamma}{p} \in (1, \infty), \quad
y^{\gamma} \in Y_p(\mu) \iff \frac{1+\alpha+\gamma}{p} \in (-\infty, \alpha).
\end{gather*}
\item The classes $X_p(\mu)$ and $Y_p(\mu)$ are dual to each other in the following sense. We have $w \in X_p(\mu) \iff w^{-\frac{1}{p-1}} \in Y_{p'}(\mu)$ and $w \in Y_p(\mu) \iff w^{-\frac{1}{p-1}} \in X_{p'}(\mu)$. In fact,
\[
[w]_{X_p(\mu)} = [w^{-\frac{1}{p-1}}]_{Y_{p'}(\mu)} \quad \text{and} \quad [w]_{Y_p(\mu)} = [w^{-\frac{1}{p-1}}]_{X_{p'}(\mu)}.
\]
\item
The additional assumptions $\alpha \in (0,\infty)$ and $\alpha \in (1,\infty)$ in Theorem \ref{thm:main1} (ii) and Theorem \ref{thm:mainYp} come from the following. We have
\begin{gather*}
A_p(\R_+, \mu) \cap X_p(\mu) \neq \emptyset \iff \alpha > 0, \quad A_p(\R_+, \mu) \cap Y_p(\mu) \neq \emptyset \iff \alpha > 0, \\
A_p(\R_+, \mu) \cap X_p(\mu) \cap Y_p(\mu) \neq \emptyset \iff \alpha > 1.
\end{gather*}
The $\Leftarrow$ implications follow by (1). The $\Rightarrow$ parts follow from H\"older's inequality. See the beginning of Section \ref{sec:lowerterms} for details.
\item
The classes $X_p(\mu)$ and $Y_p(\mu)$ are special cases of general classes of weights introduced by Muckenhoupt in \cite{Muckenhoupt1972}, in which the author obtained a complete characterization of weights admitting certain weighted Hardy's inequalities on a half-line. The class of $X_p(\mu)$ weights is characterized by the validity of the weighted Hardy's inequality $\Vert u/x_d \Vert_{L_{q,p}(\Omega_T, wd\mu)} \leq N \Vert D_d u \Vert_{L_{q,p}(\Omega_T, wd\mu)}$. See Lemma \ref{lem:Hardyinfty} for a precise statement. We use this inequality in the proofs of Theorem \ref{thm:main1} (ii) and Theorem \ref{thm:mainYp} (ii).
\item
In Theorem \ref{thm:main1} we assume $F_0=0$, whereas we do not have this restriction in Theorem \ref{thm:mainYp}. In Theorem \ref{thm:mainYp}, we rewrite $x_d^{\alpha-1}F_0=D_d(x_d^{\alpha}G_0)$ with $\lim_{x_d \to 0^+}x_d^{\alpha}G_0=0$. The class of $Y_p(\mu)$ weights is characterized by the validity of the weighted Hardy's inequality $\Vert G_0 \Vert_{L_{q,p}(\Omega_T, wd\mu)} \leq N \Vert F_0 \Vert_{L_{q,p}(\Omega_T, wd\mu)}$. See Lemma \ref{lem:Hardy0} for a precise statement. We also note that a class of weights similar to our $Y_p(\mu)$ class was used in \cite{DongPhan2021RevMatIber, DongPhanIndiana2023}, where it was denoted by $M_p(\mu)$ and used in the study of equations in non-divergence form with the conormal boundary condition. In fact, the $Y_p(\mu)$ and $M_p(\mu)$ conditions are related by $[w]_{Y_p(\mu)} = [y^p w]_{M_p(\mu)}$.
\item
In Theorem \ref{thm:mainYp}, from the $Y_p(\mu)$ condition it also follows that any weak solution $u$ to (\ref{eq:eqpar}) with $u,Du,F,f \in L_{q,p}(\Omega_T,w d\mu)$ automatically satisfies the conormal boundary condition \eqref{eq:conormal}. See the end of Section \ref{sec:lowerterms} for the proof. This fact is of independent interest, although we do not use it later in the proofs.
\end{enumerate}
\end{remark}

We now turn our attention to the corresponding elliptic systems on the half-space
 \begin{equation} \label{eq:eqell}
 \begin{gathered}
- D_i(x_d^{\alpha}A^{ij}D_j u+x_d^{\alpha-1}B^i u) + x_d^{\alpha-1}\Hat{B}^i D_i u + x_d^{\alpha-2} C u + \lambda x_d^{\alpha} C_0 u  \\
  = -D_i(x_d^{\alpha}F^i) + x_d^{\alpha-1} F_0 + \sqrt{\lambda}x_d^{\alpha} f 
 \end{gathered} \quad \text{in $\R^d_+$}
 \end{equation}
with the conormal boundary condition
\begin{equation} \label{eq:conormalell}
\lim_{x_d \to 0^+} x_d^{\alpha}(A^{dj}D_j u+ x_d^{-1}B^d u - F^d) = 0 \quad \text{on $\R^{d-1}$}.
\end{equation}
Here $A^{ij}, C_0, B^i, \Hat{B}^i, C:\R^d_+ \to \R^{n\times n}$ and $u,F^i,F_0,f:\R^d_+ \to \R^{n}$ are functions of $x \in \R^d_+$. We say that $u$ is a weak solution to (\ref{eq:eqell})-(\ref{eq:conormalell}) if
\begin{equation*}
\begin{gathered}
\int_{\R^d_+}((A^{ij}D_j u + x_d^{-1} B^i u) \cdot D_i \varphi + x_d^{-1} \Hat{B}^i D_i u \cdot \varphi + x_d^{-2} C u \cdot \varphi + \lambda C_0 u \cdot \varphi) d\mu \\
= \int_{\R^d_+}(F^i \cdot D_i \varphi + x_d^{-1} F_0 \cdot \varphi + \sqrt{\lambda} f\cdot \varphi) d\mu
\end{gathered}
\end{equation*}
holds for all $\varphi \in C_0^{\infty}(\R^d; \R^n)$. We only state the following elliptic analogue of Theorem \ref{thm:main1}. The elliptic analogue of Theorem \ref{thm:mainYp} can be stated similarly.

\begin{thm} \label{thm:ell}
Let $\alpha \in (-1, \infty), p \in (1,\infty), \kappa \in (0,1], K \in [1,\infty), R_0 \in (0,\infty]$, and $K_1 \in [0,\infty)$. Let $P=(d,\alpha,p,\kappa,K)$ denote the collection of listed parameters. There exist $\gamma_0=\gamma_0(P)>0$ and $\lambda_0=\lambda_0(P)>0$ such that the following assertions hold. Let $w = w(x)\in A_p(\R^d_+)$ with $[w]_{A_p(\R^d_+,\mu)} \leq K$. Suppose that \eqref{eq:ellipticity}, \eqref{eq:nondega0c0}, and Assumption {\rm \ref{assump}} $(\gamma_0,R_0)$ are satisfied. Also assume that $F_0=0$.

\begin{enumerate}
\item[(i)]
Suppose that $|B| \leq K_1x_d, |\Hat{B}| \leq K_1x_d$, and $|C| \leq K_1^2 x_d^2$. There exist $M=M(P)>0$ and $N=N(P)>0$ such that for any $F \in L_{p}(\R^d_+,w d\mu; \R^n)^{d}$, $f \in L_{p}(\R^d_+,w d\mu; \R^n)$, $u \in W^1_{p}(\R^d_+,w d\mu; \R^n)$, and $\lambda \geq \max(\lambda_0 R_0^{-2}, MK_1^2)$ satisfying \eqref{eq:eqell}-\eqref{eq:conormalell}, it holds that
\begin{equation*} \label{eq:estmain1ell}
\Vert Du \Vert_{L_{p}(\R^d_+,w d\mu)} + \sqrt{\lambda}\Vert u \Vert_{L_{p}(\R^d_+,w d\mu)} \leq N \Vert F\Vert_{L_{p}(\R^d_+,w d\mu)} + N \Vert f\Vert_{L_{p}(\R^d_+,w d\mu)}.
\end{equation*}
Moreover, for any $F \in L_{p}(\R^d_+,w d\mu; \R^n)^{d}$, $f \in L_{p}(\R^d_+,w d\mu; \R^n)$, and $\lambda \geq \max(\lambda_0 R_0^{-2}, MK_1^2)$ with $\lambda>0$, there exists a unique solution $u \in W^1_{p}(\R^d_+,w d\mu; \R^n)$ to \eqref{eq:eqell}-\eqref{eq:conormalell}.
\item[(ii)]
Assume, additionally, that $\alpha \in (0,\infty)$, and $w(x) = w_2(x')w_3(x_d)$, where $w_2 \in A_p(\R^{d-1}), w_3 \in A_p(\R_+, \mu) \cap X_p(\mu)$ with $[w_2]_{A_p(\R^{d-1})}\cdot [w_3]_{A_p(\R_+,\mu)} \leq K, [w_3]_{X_p(\mu)} \leq K_2, K_2 \in (0,\infty)$. There exist $\varepsilon=\varepsilon(P) > 0, M=M(P)>0$, and $N=N(P)>0$ such that the result of {\rm (i)} holds under the assumption
$$
|B| \leq \max(\varepsilon K_2^{-1}, K_1 x_d),\quad 
|\Hat{B}| \leq K_1 x_d,\quad \text{and}\quad
|C| \leq \max(K_1 K_2^{-1} x_d, K_1^2x_d^2).
$$ 
Moreover, in place of \eqref{eq:estmain1ell} holds the estimate
\begin{equation*}
\begin{aligned}
\Vert Du \Vert_{L_{p}(\R^d_+,w d\mu)} + K_2^{-1}\Vert u/x_d \Vert_{L_{p}(\R^d_+,w d\mu)} &+ \sqrt{\lambda}\Vert u \Vert_{L_{p}(\R^d_+,w d\mu)} \\
 &\leq N \Vert F\Vert_{L_{p}(\R^d_+,w d\mu)} + N \Vert f\Vert_{L_{p}(\R^d_+,w d\mu)}.
\end{aligned}
\end{equation*}
\end{enumerate}
\end{thm}

Theorem \ref{thm:ell} can be derived from Theorem \ref{thm:main1} by viewing solutions to elliptic systems as steady state solutions of the corresponding parabolic systems. More precisely, elliptic a-priori estimates above can be derived from estimates (\ref{eq:estmain1}) and (\ref{eq:estmainXp}) in Theorem \ref{thm:main1}. See, for example, the proofs of \cite[Theorem 2.6]{Krylov2007CPDE} or \cite[Theorem 1.2]{DongPhan2021RevMatIber} for details. The existence of solutions can be obtained by following the proof of existence in Theorem \ref{thm:main1}. We therefore omit the proof of Theorem \ref{thm:ell}.

\section{\texorpdfstring{$L_2$}{} solvability} \label{sec:L2}

In this section, we establish the $L_2(\Omega_T, \mu; \R^n)$ solvability for systems
\begin{equation} \label{eq:eqL2}
x_d^{\alpha}a_0(x) u_t - D_i(x_d^{\alpha}A^{ij}D_j u) + \lambda x_d^{\alpha} C_0 u 
  = - D_i(x_d^{\alpha}F^i) + \sqrt{\lambda}x_d^{\alpha} f \quad \text{in } \Omega_T
\end{equation}
with the conormal boundary condition
\begin{equation} \label{eq:L2conorm}
\lim_{x_d \to 0^+} x_d^{\alpha}(A^{dj}D_j u - F^d) = 0 \quad \text{on } (-\infty, T) \times \R^{d-1}.
\end{equation}
In this section, the coefficient $a_0=a_0(x)$ may depend on all the $x$ variables, and we do not impose any regularity assumptions on the coefficients $a_0(x), A^{ij}(t,x)$, and $C_0(t,x)$ other than (\ref{eq:ellipticity}) and (\ref{eq:nondega0c0}). The following theorem is the main result of this section. Denote $\R_T = (-\infty, T)$.

\begin{thm}[$L_2$ solvability] \label{thm:L2solv}
Let $\alpha \in (-1,\infty), \kappa \in (0, 1], \lambda \geq 0$, and suppose that \eqref{eq:ellipticity} and \eqref{eq:nondega0c0} are satisfied. If $u \in \mathcal{H}^{1}_2(\Omega_T, \mu; \R^n)$ is a weak solution to \eqref{eq:eqL2}-\eqref{eq:L2conorm} for some $F \in L_2(\Omega_T,\mu; \R^n)^{d}$ and $f \in L_2(\Omega_T,\mu; \R^n)$, then $u \in C_b(\overline{\R_T}; L_2(\R^d_+, \mu; \R^n))$ and
\begin{equation} \label{eq:L2estimate}
\sup_{t \in \overline{\R_T}}\int_{\R^d_+} |u(t,\cdot)|^2 d\mu + \int_{\Omega_T} |Du|^2 d\mu + \lambda \int_{\Omega_T} |u|^2 d\mu \leq N \int_{\Omega_T} |F|^2 d\mu + N \int_{\Omega_T} |f|^2 d\mu
\end{equation}
for some $N=N(\kappa)>0$. Moreover, for any $\lambda > 0$, $F \in L_2(\Omega_T,\mu; \R^n)^{d}$, and $f \in L_2(\Omega_T,\mu; \R^n)$, there exists a unique weak solution $u \in \mathcal{H}^{1}_2(\Omega_T,\mu; \R^n)$ to \eqref{eq:eqL2}-\eqref{eq:L2conorm}.
\end{thm}

Here $C_b(\overline{\R_T}; L_2(\R^d_+, \mu; \R^n))$ stands for the space of continuous $L_2(\R^d_+, \mu)$-valued functions on $\overline{\R_T}$. Before giving the proof, we introduce some notation and auxiliary results. 
We first recall the Steklov averages in the time variable. For $h>0$, the backward Steklov average of a function $u(t,x)$ is the function
\[
u_h(t,x) := \frac{1}{h}\int_{t-h}^t u(s,x)ds = \int_0^1 u(t-sh,x)ds.
\]
We will use the following properties of Steklov averages. If $u \in L_2(\Omega_T, \mu; \R^n)$, then
\begin{gather*}
\Vert u_h \Vert_{L_2(\Omega_T, \mu; \R^n)} \leq \Vert u \Vert_{L_2(\Omega_T, \mu; \R^n)} \quad \text{for any } h>0, \\
(\partial_t u_h)(t,x) = \frac{u(t,x)-u(t-h,x)}{h} \quad \text{for any } h>0, \\
\lim_{h \to 0^+} \Vert u_h - u \Vert_{L_2(\Omega_T, \mu; \R^n)} = 0.
\end{gather*}

For $T = \infty$, let $\Hat{u}(\tau,x) := (2\pi)^{-1/2} \int_{\R} u(t,x) e^{-i t \tau} dt$ be the Fourier transform in the time variable of $u \in L_2(\R^{d+1}_+, \mu; \C^n) \cap L_1(\R;L_2(\R^d_+,\mu; \C^n))$. By density and Plancherel's theorem, the Fourier transform in time extends to an isometric isomorphism $\mathcal{F}_t:L_2(\R^{d+1}_+, \mu; \C^n) \to L_2(\R^{d+1}_+, \mu; \C^n)$. The Hilbert transform in time $H$ is the isometric isomorphism $L_2(\R^{d+1}_+, \mu; \C^n) \to L_2(\R^{d+1}_+, \mu; \C^n)$ defined by
\begin{equation*}
\widehat{Hu}(\tau,x) = -i \operatorname{sgn}\tau \; \Hat{u}(\tau, x).
\end{equation*}
For functions $u \in L_2(\R^{d+1}_+, \mu; \C^n)$ with $|\tau|^{1/2} \Hat{u} \in L_2(\R^{d+1}_+, \mu; \C^n)$, the $(-\partial_t^2)^{1/4}$ operator is defined by
\[
(-\partial_t^2)^{1/4}u = \mathcal{F}_t^{-1}(|\tau|^{1/2} \Hat{u}) \in L_2(\R^{d+1}_+, \mu; \C^n).
\]

A straightforward verification shows that if $u \in L_2(\R^{d+1}_+, \mu; \R^n)$ is $\R^n$-valued, then $Hu \in L_2(\R^{d+1}_+, \mu; \R^n)$ is also $\R^n$-valued. Moreover, if $u \in L_2(\R^{d+1}_+, \mu; \R^n)$ is $\R^n$-valued and $|\tau|^{1/2} \Hat{u} \in L_2(\R^{d+1}_+, \mu)$, then $(-\partial_t^2)^{1/4}u \in L_2(\R^{d+1}_+, \mu; \R^n)$ is also $\R^n$-valued.

To prove the existence of solutions in Theorem \ref{thm:L2solv}, we consider the following slightly more general systems with half-time derivatives in the spirit of \cite{KaplanPisa1966} and \cite{JungKim2025JDE}:
\begin{equation} \label{eq:1/2time}
x_d^{\alpha}a_0(x) u_t - D_i(x_d^{\alpha}A^{ij}D_j u) + \lambda x_d^{\alpha} C_0 u 
  = x_d^{\alpha} (-\partial_t^2)^{1/4}g - D_i(x_d^{\alpha}F^i) + \sqrt{\lambda}x_d^{\alpha} f \quad \text{in } \R^{d+1}_+.
\end{equation}
The solution space to (\ref{eq:1/2time})-(\ref{eq:L2conorm}) is the (real) Hilbert space
\[
H^{1/2,1}_2(\R^{d+1}_+, \mu; \R^n) = \big\{u \in L_2(\R^{d+1}_+, \mu; \R^n): (-\partial_t^2)^{1/4}u, D_1u,\dots,D_d u \in L_2(\R^{d+1}_+, \mu; \R^n)\big\}
\]
with the inner product given by
\[
\langle u, v \rangle_{H^{1/2,1}_2(\R^{d+1}_+, \mu; \R^n)} = \int_{\R^{d+1}_+}\left[(-\partial_t^2)^{1/4}u \cdot (-\partial_t^2)^{1/4}v + D_iu \cdot D_i v + u \cdot v\right]d\mu.
\]
We say that $u \in H^{1/2,1}_2(\R^{d+1}_+, \mu; \R^n)$ is an $H^{1/2,1}_2(\R^{d+1}_+, \mu; \R^n)$-solution to (\ref{eq:1/2time})-(\ref{eq:L2conorm}) if 
\begin{gather*}
\int_{\R^{d+1}_+}\Bigl(-a_0(x) [H(-\partial_t^2)^{1/4}u] \cdot (-\partial_t^2)^{1/4}v + A^{ij}D_j u \cdot D_i v + \lambda C_0 u \cdot v \Bigr) d\mu \\
= \int_{\R^{d+1}_+}\Bigl(g \cdot (-\partial_t^2)^{1/4} v + F^i \cdot D_i v + \lambda^{1/2} f \cdot v  \Bigr) d\mu
\end{gather*}
holds for all $v \in H^{1/2,1}_2(\R^{d+1}_+, \mu; \R^n)$.
Clearly, if $g=0$, an $H^{1/2,1}_2(\R^{d+1}_+, \mu; \R^n)$-solution is also a weak solution to (\ref{eq:eqL2})-(\ref{eq:L2conorm}) in the sense of identity (\ref{eq:soldefvar}). We will need the following result on the solvability of \eqref{eq:1/2time}-\eqref{eq:L2conorm} in $H^{1/2,1}_2(\R^{d+1}_+, \mu; \R^n)$.

\begin{thm} \label{thm:1/2timeL2}
Let $\alpha \in (-1,\infty)$ and suppose that \eqref{eq:ellipticity} and \eqref{eq:nondega0c0} are satisfied. Then for any $\lambda > 0, F \in L_2(\R^{d+1}_+, \mu; \R^n)^d,$ and $f, g \in L_2(\R^{d+1}_+, \mu; \R^n)$, there exists a unique $H^{1/2,1}_2(\R^{d+1}_+, \mu; \R^n)$-solution $u$ to \eqref{eq:1/2time}-\eqref{eq:L2conorm}.
\end{thm}

\begin{proof}
The proof uses a parabolic version of the Lax-Milgram theorem. We only note that the symmetry condition $a_0=a_0^T$ ensures the validity of the identity
\[
\int_{\R^{d+1}_+}a_0(x)Hf \cdot f d\mu = 0
\]
for all $f \in L_2(\R^{d+1}_+,\mu; \R^n)$. For the rest of the details, we refer to the proof of \cite[Theorem 4.1]{JungKim2025JDE}.
\end{proof}

\begin{proof}[Proof of Theorem {\rm \ref{thm:L2solv}}]
By taking Steklov averages, for any $h>0$ and $\varphi \in C_0^{\infty}((-\infty, T) \times \R^d; \R^n)$, 
\[
\int_{\Omega_T} [-a_0(x)u_h \cdot \varphi_t + (A^{ij}D_j u)_h \cdot D_i \varphi + \lambda(C_0 u)_h \cdot \varphi]d\mu = \int_{\Omega_T} [F^i_h \cdot D_i \varphi + \lambda^{1/2} f_h \cdot \varphi] d\mu.
\]
Here we can integrate by parts in $t$ to get
\[
\int_{\Omega_T} [a_0(x)\partial_t u_h \cdot \varphi + (A^{ij}D_j u)_h \cdot D_i \varphi + \lambda(C_0 u)_h \cdot \varphi]d\mu = \int_{\Omega_T} [F^i_h \cdot D_i \varphi + \lambda^{1/2} f_h \cdot \varphi] d\mu.
\]
By considering cut-off in time, for any $v \in C_0^{\infty}(\R^{d+1}; \R^n)$ and $S\leq T$,
\begin{equation} \label{eq:1}
\int_{\Omega_S} [a_0(x)\partial_t u_h \cdot v + (A^{ij}D_j u)_h \cdot D_i v + \lambda(C_0 u)_h \cdot v]d\mu = \int_{\Omega_S} [F^i_h \cdot D_i v + \lambda^{1/2} f_h \cdot v] d\mu.
\end{equation}
By approximation, (\ref{eq:1}) holds for all $v \in L_2(\Omega_S, \mu; \R^n)$ with $Dv \in L_2(\Omega_S, \mu)$.

Now take $v = u_h$. Then for $S \in \overline{\R_T}$ holds
\begin{equation} \label{eq:3}
\begin{aligned}
&\int_{\Omega_S} a_0(x) \partial_t u_h \cdot u_h d\mu = [\text{using $a_0 = a_0^T$}] = \int_{\Omega_S}\frac{1}{2} \partial_t (a_0(x)u_h \cdot u_h) d\mu \\
&= \int_{\R^d_+} \frac{1}{2}a_0(x) u_h(S,\cdot) \cdot u_h(S,\cdot) d\mu \geq [\text{by (\ref{eq:nondega0c0})}] \geq \frac{\kappa}{2}\int_{\R^d_+} |u_h|^2(S,\cdot) d\mu \geq 0.
\end{aligned}
\end{equation}
If $S=T=\infty$, then
\[
\int_{\R^{d+1}_+} a_0(x) \partial_t u_h \cdot u_h d\mu = 0.
\]

Therefore, for all $S \leq T$ holds
\[
\int_{\Omega_S} [(A^{ij}D_j u)_h \cdot D_i u_h + \lambda(C_0 u)_h \cdot u_h]d\mu \leq \int_{\Omega_S} [F^i_h \cdot D_i u_h + \lambda^{1/2} f_h \cdot u_h] d\mu.
\]
In the limit $h \to 0^+$, we get
\[
\int_{\Omega_S} [A^{ij}D_j u \cdot D_i u + \lambda C_0 u \cdot u]d\mu \leq \int_{\Omega_S} [F^i \cdot D_i u + \lambda^{1/2} f \cdot u] d\mu.
\]
By (\ref{eq:ellipticity}), (\ref{eq:nondega0c0}), and Young's inequality, we get
\begin{equation} \label{eq:4}
\int_{\Omega_S} |Du|^2 d\mu + \lambda \int_{\Omega_S} |u|^2 d\mu \leq N \int_{\Omega_S} |F|^2 d\mu + N \int_{\Omega_S} |f|^2 d\mu.
\end{equation}
Here and throughout the proof $N=N(\kappa)>0$. From (\ref{eq:1}), (\ref{eq:3}), and (\ref{eq:4}), for all $S \in \overline{\R_T}$ it holds that ($\Vert \cdot \Vert = \Vert \cdot \Vert_{L_2(\Omega_S,\mu)}$)
\begin{equation} \label{eq:2}
\begin{aligned}
&\frac{\kappa}{2}\Vert u_h(S,\cdot) \Vert_{L_2(\R^d_+,\mu)} \leq \int_{\Omega_S} a_0(x) \partial_t u_h \cdot u_h d\mu \\
&\leq (\Vert F^i_h \Vert + \Vert (A^{ij}D_j u)_h \Vert) \Vert D_i u_h \Vert
+ \lambda^{1/2} \Vert u_h \Vert (\Vert f_h \Vert + \lambda^{1/2}\Vert (C_0 u)_h \Vert) \\
&\leq N \Vert Du \Vert^2 + N \lambda \Vert u \Vert^2 + \Vert F \Vert^2 + \Vert f \Vert^2 
\leq N \Vert F \Vert_{L_2(\Omega_S,\mu)}^2 + N \Vert f \Vert_{L_2(\Omega_S,\mu)}^2.
\end{aligned}
\end{equation}
From (\ref{eq:1}), for all $h_1,h_2 > 0$, we have
\begin{align*}
&\frac{1}{2}\int_{\R^d_+}a_0(x)(u_{h_1}-u_{h_2})\cdot(u_{h_1}-u_{h_2})(S,x)d\mu(x) \\
&= \int_{\Omega_S} \left[ \Bigl(F^i_{h_1}-F^i_{h_2}-(A^{ij}D_j u)_{h_1}+(A^{ij}D_j u)_{h_2}\Bigr)\cdot(D_i u_{h_1}-D_i u_{h_2}) \right.\\
&\left.\quad + \lambda^{1/2}\Bigl(f_{h_1}-f_{h_2}-\lambda^{1/2}(C_0 u)_{h_1}+\lambda^{1/2}(C_0 u)_{h_2}\Bigr)\cdot(u_{h_1}-u_{h_2})\right]d\mu.
\end{align*}
Hence $\Vert u_{h_1} - u_{h_2} \Vert_{C_b(\overline{\R_T}; L_2(\R^d_+, \mu))} \to 0$ as $h_1,h_2 \to 0^+$, from which and also (\ref{eq:2}), we have
\[
u \in C_b(\overline{\R_T}; L_2(\R^d_+, \mu))
\]
and
\[
\sup_{t \in \overline{\R_T}} \int_{\R_d^+}|u(t,\cdot)|^2 d\mu \leq N\Vert F \Vert_{L_2(\Omega_T,\mu)}^2 + N \Vert f \Vert_{L_2(\Omega_T,\mu)}^2.
\]
Estimate (\ref{eq:L2estimate}) is proved. By extension, it suffices to prove the existence of solutions for the case $T=\infty$, which follows by Theorem \ref{thm:1/2timeL2}.
\end{proof}

\section{Systems with simple coefficients} \label{sec:4}
In this section, we study systems with simple coefficients. Throughout this section, we assume that the coefficients $a_0=a_0(x_d), A^{ij}=A^{ij}(x_d)$, and $C_0=C_0(x_d)$ depend on $x_d$ only, $B = \Hat{B} = 0$, $C = 0$, and $F_0 = 0$. That is, we consider systems
\begin{equation} \label{eq:simplenonhom}
x_d^{\alpha} a_0(x_d) u_t - D_i(x_d^{\alpha}A^{ij}(x_d)D_j u) + \lambda x_d^{\alpha} C_0(x_d) u = -D_i(x_d^{\alpha} F^i) +\sqrt{\lambda}x_d^{\alpha} f
\end{equation}
with the conormal boundary condition
\begin{equation} \label{eq:conormalsimplenonhom}
\lim_{x_d \to 0^+} x_d^{\alpha}(A^{dj}(x_d)D_j u - F^d) = 0.
\end{equation}
\noindent
The goal of this section is to prove the following theorem, which is a special case of Theorem \ref{thm:main1} (i).

\begin{thm} \label{thm:exisuniqsimpleunweighted}
Let $\alpha \in (-1,\infty), p \in (1,\infty), \kappa \in (0,1]$, and $\lambda>0$. Suppose that \eqref{eq:ellipticity} and \eqref{eq:nondega0c0} are satisfied. Then, for any $F \in L_p(\Omega_T,\mu; \R^n)^d$ and $f \in L_p(\Omega_T,\mu; \R^n)$, there exists a unique solution $u \in \mathcal{H}^1_p(\Omega_T,\mu; \R^n)$ to \eqref{eq:simplenonhom}-\eqref{eq:conormalsimplenonhom}. Moreover,
\begin{equation} \label{eq:apriorisimpleunweighted}
\Vert Du \Vert_{L_p(\Omega_T,\mu)} + \sqrt{\lambda}\Vert u \Vert_{L_p(\Omega_T,\mu)} \leq N \Vert F\Vert_{L_p(\Omega_T,\mu)} + N \Vert f\Vert_{L_p(\Omega_T,\mu)},
\end{equation}
where $N=N(d,\alpha,p,\kappa)>0$.
\end{thm}

Throughout this section, it is assumed that $\alpha \in (-1,\infty)$, $\lambda \geq 0$, and that \eqref{eq:ellipticity} and \eqref{eq:nondega0c0} are satisfied with some $\kappa \in (0,1]$.
For compactness of notation, denote
\[
\mathcal{P}_0u = x_d^{\alpha} a_0(x_d) u_t - D_i(x_d^{\alpha}A^{ij}(x_d)D_j u) + \lambda x_d^{\alpha} C_0(x_d) u.
\]
\subsection{Lipschitz estimates for homogeneous systems}

Let $z_0 =(t_0, x_0) \in \overline{\R^{d+1}_+}$ and $R > 0$. We consider the homogeneous system
\begin{equation} \label{eq:simplehom}
x_d^{\alpha} a_0(x_d) u_t - D_i(x_d^{\alpha}A^{ij}(x_d)D_j u) + \lambda x_d^{\alpha} C_0(x_d) u = 0 \quad \text{in } \; Q_R^+(z_0),
\end{equation}
with the conormal boundary condition
\begin{equation} \label{eq:conormalsimplehom}
\lim_{x_d \to 0^+} x_d^{\alpha} A^{dj}(x_d)D_j u = 0 \quad \text{on $Q_R'(z_0')$ if $x_{0d} \leq R$}.
\end{equation}

Our goal is to derive Lipschitz estimates for (\ref{eq:simplehom})-(\ref{eq:conormalsimplehom}). We begin with the following Caccioppoli-type lemma.

\begin{lem} \label{lem:caccioppoli}
Let $0 < r < R, z_0 = (t_0,x_0) \in \overline{\R^{d+1}_+}$, and $u \in \mathcal{H}^1_2(Q_R^+(z_0), \mu; \R^n)$ be a weak solution to \eqref{eq:simplehom}-\eqref{eq:conormalsimplehom} in $Q_R^+(z_0)$. Then $u_t \in L_2(Q_r^+(z_0))$ and
\begin{equation*}
\begin{gathered}
\int_{Q_r^+(z_0)} (|Du|^2+\lambda |u|^2) d\mu \leq \frac{N}{(R-r)^2} \int_{Q_R^+(z_0)} |u|^2 d\mu, \\
\int_{Q_r^+(z_0)} |u_t|^2 d\mu \leq \frac{N}{(R-r)^2} \int_{Q_R^+(z_0)} (|Du|^2+\lambda |u|^2) d\mu,
\end{gathered}
\end{equation*}
where $N=N(\kappa) > 0$.
\end{lem}
The proof is standard; see, for example, the argument in \cite[p. 1082-1084]{DongKimSIAM2011}. For the first inequality, we test the equation by $u\varphi^2$, where $\varphi$ is an appropriately chosen smooth cutoff function. As in the proof of Theorem \ref{thm:L2solv}, one can use the Steklov averages to justify taking $u\varphi^2$ as a test function. For the second inequality, note that Steklov averages $u_h$ are solutions to \eqref{eq:simplehom}-\eqref{eq:conormalsimplehom} on $(t_0-R^2+h,t_0)\times D_{R}^+(x_0)$. We then test the equation for $u_h$ by $(\partial_t u_h) \varphi_k^2$, where $\varphi_k$ are appropriately chosen smooth cutoff functions with expanding supports, and use an iteration argument combined with the first inequality. This proves the second inequality with $u_h$ instead of $u$ on the left hand side. We then pass to a limit along a subsequence as $h \to 0^+$. Note that the second inequality holds if the coefficients are merely measurable and independent of $t$, and for the first inequality, we used only that $a_0$ does not depend on $t$.

When the coefficients are independent of $t$ and $x'$, we can iterate the estimates in Lemma \ref{lem:caccioppoli} to show that all partial derivatives of $u$ in $t,x'$ exist and are again solutions to \eqref{eq:simplehom}-\eqref{eq:conormalsimplehom}. For $I,J \geq 0$, we denote
\[
|\partial_t^I D_{x'}^J u|^2 = \sum_{j_1,\dots,j_J =1}^{d-1} |\partial_t^I D_{j_1,\dots,j_J}u|^2, \quad |D\partial_t^I D_{x'}^J u|^2 = \sum_{j=1}^d \sum_{j_1,\dots,j_J = 1}^{d-1} |\partial_t^I D_{j,j_1,\dots,j_J}u|^2.
\]

For $\varepsilon > 0$, the mollification of a function $u(t,x)$ in the $t,x'$ variables is the function
\[
u^{(\varepsilon)}(t,x',x_d) = \int_0^1 \int_{B_1'} u(t-\varepsilon^2 s, x'-\varepsilon y',x_d) \eta_1(s) \eta_2(y') dy'ds,
\]
where $\eta_1 \in C_0^{\infty}(0,1), \eta_2 \in C_0^{\infty}(B_1')$ are fixed non-negative functions with unit integrals.

\begin{lem} \label{lem:caccioimproved}
In the setting of Lemma {\rm \ref{lem:caccioppoli}}, for any $I,J \geq 0$, we have that $\partial_t^I D_{x'}^J u \in \mathcal{H}^1_2(Q_r^+(z_0), \mu)$ and $\partial_t^I D_{x'}^J u$ is a weak solution to \eqref{eq:simplehom}-\eqref{eq:conormalsimplehom} in $Q_R^+(z_0)$. Moreover,
\begin{gather*} 
\int_{Q_r^+(z_0)} |\partial_t^I D_{x'}^J u|^2 d\mu \leq \frac{N_{\kappa,I,J}}{(R-r)^{4I+2J-2}} \int_{Q_R^+(z_0)} (|Du|^2+\lambda |u|^2) d\mu, \quad I+J \geq 1, \\
\int_{Q_r^+(z_0)} (\lambda^k |D\partial_t^I D_{x'}^J u|^2+\lambda^{k+1} |\partial_t^I D_{x'}^J u|^2) d\mu \leq \frac{N_{\kappa,I,J,k}}{(R-r)^{4I+2J+2k}} \int_{Q_R^+(z_0)} (|Du|^2+\lambda |u|^2) d\mu, \, I,J,k \geq 0.
\end{gather*}
\end{lem}

These inequalities are first proved with $u^{(\varepsilon)}$ instead of $u$ in the left hand sides for all small enough $\varepsilon$ by iterating the estimates in Lemma \ref{lem:caccioppoli}. We then pass to a limit along a subsequence as $\varepsilon \to 0^+$.

Before we proceed further, we introduce some notation and auxiliary inequalities. Denote
\[
\mathcal{U} = A^{dj} D_j u.
\]
When estimating $Du$, from the structure of the equations, it is often convenient to first estimate $\mathcal{U}$ and $D_{x'}u$ separately. From (\ref{eq:ellipticity}), we have
\begin{equation} \label{eq:estDduCalU}
\begin{gathered}
\kappa|D_d u| \leq |A^{dd}D_d u|= |\mathcal{U}-\sum_{j=1}^{d-1}A^{dj}D_j u| \leq |\mathcal{U}|+\kappa^{-1}|D_{x'}u|, \\
|\mathcal{U}| \leq \kappa^{-1}|Du|.
\end{gathered}
\end{equation}

For $p \in [1,\infty)$, even $k \geq 0$, $r>0$, and $u:Q_r' \to \R^n$, we set
\[
\Vert u \Vert^p_{W^{k/2,k}_p(Q_r')} := \int_{Q_r'} \sum_{2I+J \leq k} |\partial_t^I D_{x'}^J u|^p dx'dt.
\]
When $k > (d+1)/p$, we have the Sobolev embedding $W^{k/2,k}_p(Q_r';\R^n) \hookrightarrow L_{\infty}(Q_r'; \R^n)$ with
\[
\sup_{Q_r'}|u|^p \leq N_{d,p} \fint_{Q_r'}\sum_{2I+J\leq k}|r^{2I+J} \partial_t^I D_{x'}^J u|^p dx'dt.
\]
For scalar-valued functions the above inequality is well-known. For $\R^n$-valued functions $u$, one can apply the inequality to scalar functions $u(t,x') \cdot \xi$, where $\xi \in \R^n$ is fixed, and take $\sup$ over all $\xi$ with $|\xi|=1$. Hence, the constant in the above Sobolev inequality does not depend on $n$.

Lastly, denote $d\mu_{\alpha} = s^{\alpha}ds$ on $\R_+$. A direct verification shows that whenever $0 \leq a_1 \leq a_2 < b_2 \leq b_1$, we have
\begin{equation} \label{eq:mudoubling1d}
\frac{\mu_{\alpha}((a_1,b_1))}{\mu_{\alpha}((a_2,b_2))} = \frac{b_1^{1+\alpha}-a_1^{1+\alpha}}{b_2^{1+\alpha}-a_2^{1+\alpha}} \leq N_{\alpha} \left(\frac{b_1-a_1}{b_2-a_2}\right)^{1+\alpha_+},
\end{equation}
where $N_{\alpha} = \max((1+\alpha)^{-1},1)$ is sharp. In particular, the measure $\mu_{\alpha}$ is doubling.

In the following lemma, we prove interior reverse H\"older type $L_{\infty}$ estimates of solutions.
\begin{lem} \label{lem:intbddnssest}
Let $p \in (1,2], r>0$, $z_0 \in \R^{d+1}_+$ with $x_{0d} \geq 2r$, and let $u \in \mathcal{H}^1_p(Q_r(z_0), \mu; \R^n)$ be a weak solution to \eqref{eq:simplehom} in $Q_r(z_0)$. Then for any $\tau \in (0,1)$ it holds that
\begin{equation} \label{eq:intbddnessestsimplecoeffhom}
\begin{gathered}
\Vert u \Vert_{L_{\infty}(Q_{\tau r}(z_0))} \leq N (1-\tau)^{-\frac{d+2}{p}}\left( \fint_{Q_r(z_0)} |u|^p d\mu \right)^{1/p}, \\
\Vert Du \Vert_{L_{\infty}(Q_{\tau r}(z_0))} + \sqrt{\lambda}\Vert u \Vert_{L_{\infty}(Q_{\tau r}(z_0))} \leq N (1-\tau)^{-\frac{d+2}{p}}\left( \fint_{Q_r(z_0)} (|Du|^p+\lambda^{p/2} |u|^p) d\mu \right)^{1/p},
\end{gathered}
\end{equation}
where $N=N(d,\alpha,p,\kappa) > 0$.
\end{lem}

\begin{proof} Since $1/2 \leq x_{d}/x_{0d} \leq 3/2$ in $Q_r(z_0)$, by multiplying the equation with $x_{0d}^{-\alpha}$ and a translation of coordinates, we may assume without loss that $\alpha = 0$ and $z_0 = 0$. Moreover, it suffices to consider the case when $\tau \in [1/2,1)$.

We first consider the case $p=2$ and $\tau = 1/2$. By scaling, it is sufficient to consider $r=1$. For compactness of notation, let
\begin{equation*}
C_1 = \left( \int_{Q_1} |u(t,x)|^2 dz \right)^{1/2}, \quad
C_2 = \left( \int_{Q_1} (|Du(t,x)|^2+\lambda |u(t,x)|^2) dz \right)^{1/2}.
\end{equation*}
Fix the smallest even integer $k$ with $k > (d+1)/2$. By the Sobolev embedding and Lemmas \ref{lem:caccioppoli} and \ref{lem:caccioimproved}, for $z \in Q_{1/2}$, it holds that
\begin{equation*}
\begin{aligned}
|u(z',x_d)|^2 &\leq N \int_{-1/2}^{1/2} (|u|^2+|D_d u|^2)(z',s) ds \\
 &\leq N_d \int_{-1/2}^{1/2}(\Vert u(\cdot,s)\Vert^2_{W^{k/2,k}_2(Q_{1/2}')}+\Vert D_d u(\cdot,s)\Vert^2_{W^{k/2,k}_2(Q_{1/2}')}) ds \leq N_{d,\kappa} C_1^2.
\end{aligned}
\end{equation*}
By applying the above estimate to $D_{x'}u$ and $\sqrt{\lambda} u$ and using Lemmas \ref{lem:caccioppoli} and \ref{lem:caccioimproved}, on $Q_{1/2}$ we have
\[
|D_{x'} u| + \sqrt{\lambda} |u| \leq N_{d,\kappa} C_2.
\]
From the definition of $\mathcal{U}$ and the equation, we have $|\partial_t^I D_{x'}^J \mathcal{U}| \leq N_{\kappa}|D\partial_t^I D_{x'}^J u|$ and $|D_d \partial_t^I D_{x'}^J \mathcal{U}| \leq N_{d,\kappa} (|\partial_t \partial_t^I D_{x'}^J u|+|DD_{x'} \partial_t^I D_{x'}^J u|+\lambda|\partial_t^I D_{x'}^J u|)$ for all $I,J \geq 0$. We then estimate the term $\mathcal{U}$ similarly as in the estimate of $u$ above, and get
\[
|\mathcal{U}| \leq N_{d,\kappa} C_2.
\]
Combining the above estimates and (\ref{eq:estDduCalU}), we estimate $D_d u$ and arrive to (\ref{eq:intbddnessestsimplecoeffhom}) for $p = 2$ and $\tau = 1/2$.

By a covering argument, see, for example, \cite[p. 184-186]{GiaquintaMartinazzi} or the proof of Lemma \ref{lem:bdrbddnssest} below, estimates (\ref{eq:intbddnessestsimplecoeffhom}) hold for $p=2$ and general $\tau \in [1/2,1)$. Then by a standard iteration argument, see again \cite[p. 184-186]{GiaquintaMartinazzi}, estimates (\ref{eq:intbddnessestsimplecoeffhom}) hold for general $p \in (1,2]$, assuming that initially $u \in \mathcal{H}^1_2(Q_r)$.

Assume now that $u \in \mathcal{H}^1_p(Q_r)$, where $p \in (1,2)$, and let $\tau \in [1/2,1)$ be fixed. Let $u^{(\varepsilon)}$ be the mollification of $u$ in the $t,x'$ variables. Let $r_1 = \frac{1+\tau}{2} r$ and $r_2 = \frac{3+\tau}{4} r$. As above,
\[
\Vert u^{(\varepsilon)} \Vert_{L_{\infty}(Q_{r_1})}^p \lesssim \int_{-r_2}^{r_2} (\Vert D_d u^{(\varepsilon)}(\cdot, s) \Vert^p_{W^{k/2,k}_p(Q'_{r_2})}+\Vert u^{(\varepsilon)}(\cdot, s) \Vert^p_{W^{k/2,k}_p(Q'_{r_2})})ds < \infty
\]
for a fixed large enough even $k=k(p,d)$. Then $\Vert D_{x'}u^{(\varepsilon)} \Vert_{L_{\infty}(Q_{r_1})} < \infty$ as well. Next,
\begin{align*}
\Vert \mathcal{U}^{(\varepsilon)} \Vert^p_{L_{\infty}(Q_{r_1})} &\lesssim \int_{-r_2}^{r_2} (\Vert D_d \mathcal{U}^{(\varepsilon)}(\cdot, s) \Vert^p_{W^{k/2,k}_p(Q'_{r_2})}+\Vert \mathcal{U}^{(\varepsilon)}(\cdot, s) \Vert^p_{W^{k/2,k}_p(Q'_{r_2})})ds \\
&\lesssim \int_{-r_2}^{r_2} (\Vert D_d u^{(\varepsilon)}(\cdot, s) \Vert^p_{W^{k/2+1,k+2}_p(Q'_{r_2})}+\Vert u^{(\varepsilon)}(\cdot, s) \Vert^p_{W^{k/2+1,k+2}_p(Q'_{r_2})})ds < \infty,
\end{align*}
and therefore $\Vert D_d u^{(\varepsilon)} \Vert_{L_{\infty}(Q_{r_1})} < \infty$ as well.
In particular, $u^{(\varepsilon)} \in \mathcal{H}^1_2(Q_{r_1})$ for all small enough $\varepsilon$. Then, (\ref{eq:intbddnessestsimplecoeffhom}) holds with $u^{(\varepsilon)}$ in place of $u$ in the left hand sides. Passing to the limit along a subsequence as $\varepsilon \to 0^+$, we obtain (\ref{eq:intbddnessestsimplecoeffhom}) for $u$.
\end{proof}

The following lemma is the key intermediate step for the proofs of our main results. Here we prove boundary reverse H\"older type $L_{\infty}$ estimates of solutions.
\begin{lem} \label{lem:bdrbddnssest}
Let $p \in (1,2]$, $r>0$, and $u \in \mathcal{H}^1_p(Q_r^+, \mu; \R^n)$ be a weak solution to \eqref{eq:simplehom}-\eqref{eq:conormalsimplehom} in $Q_r^+$. Then for any $\tau \in (0,1)$ it holds that

\begin{equation} \label{eq:bdrybddnessestsimplehom}
\begin{gathered}
\Vert u \Vert_{L_{\infty}(Q_{\tau r}^+)} \leq N (1-\tau)^{-\frac{d+2+\alpha_+}{p}} \left( \fint_{Q_r^+} |u|^p d\mu \right)^{1/p}, \\
\Vert Du \Vert_{L_{\infty}(Q_{\tau r}^+)} + \sqrt{\lambda}\Vert u \Vert_{L_{\infty}(Q_{\tau r}^+)} \leq N  (1-\tau)^{-\frac{d+2+\alpha_+}{p}} \left( \fint_{Q_r^+} (|Du|^p+\lambda^{p/2} |u|^p) d\mu \right)^{1/p}, \\
\Vert \mathcal{U}/x_d \Vert_{L_{\infty}(Q_{\tau r}^+)} \leq Nr^{-1}  (1-\tau)^{-\frac{d+2+\alpha_+}{p}-1} \left( \fint_{Q_r^+} (|Du|^p+\lambda^{p/2} |u|^p) d\mu \right)^{1/p},
\end{gathered}
\end{equation}
where $N=N(d,\alpha,p,\kappa) > 0$.
\end{lem}

\begin{proof} It suffices to consider the case when $\tau \in [1/2,1)$.

We first consider the case $p=2$ and $\tau = 1/2$. By scaling, it is sufficient to consider $r=1$. For compactness of notation, let
\begin{equation*}
C_1 = \left( \int_{Q_1^+} |u(t,x)|^2 d\mu(z) \right)^{1/2}, \quad
C_2 = \left( \int_{Q_1^+} (|Du(t,x)|^2+\lambda |u(t,x)|^2) d\mu(z) \right)^{1/2}.
\end{equation*}
Throughout the proof, the constants $N$ may depend on $d,\alpha,p$, and $\kappa$. Fix the smallest even integer $k$ with $k > (d+1)/2$. From the equation, we have $|D_d (x_d^{\alpha}\mathcal{U})| \leq N x_d^{\alpha}(|u_t|+|DD_{x'}u|+\lambda|u|)$.  Then from $\lim_{x_d \to 0^+} x_d^{\alpha}\mathcal{U} = 0$, Sobolev embedding, H\"older's inequality, and Lemma \ref{lem:caccioimproved}, for $z \in Q_{1/2}^+$ it holds that
\begin{equation} \label{eq:estU1}
\begin{aligned}
&|x_d^{\alpha} \mathcal{U}(z',x_d)| \leq N \int_0^{x_d} s^{ \alpha}(|u_t|+|D_{x'}D u| + \lambda |u|)(z',s) ds \\
&\leq N \int_0^{x_d} s^{ \alpha}(\Vert u_t(\cdot,s)\Vert_{W_2^{k/2,k}(Q_{1/2}')}+\Vert D_{x'}D u(\cdot,s)\Vert_{W_2^{k/2,k}(Q_{1/2}')}+\lambda\Vert u(\cdot,s)\Vert_{W^{k/2,k}_2(Q_{1/2}')}) ds \\
&\leq N x_d^{\frac{1+\alpha}{2}}\left(\int_0^{1/2} \left[\Vert u_t(\cdot,s)\Vert^2_{W_2^{k/2,k}(Q_{1/2}')}+\Vert D_{x'}D u(\cdot,s)\Vert^2_{W_2^{k/2,k}(Q_{1/2}')}+\lambda^2\Vert u(\cdot,s)\Vert^2_{W_2^{k/2,k}(Q_{1/2}')}\right] s^{ \alpha} ds\right)^{1/2} \\
&\leq N x_d^{\frac{1+\alpha}{2}} C_2.
\end{aligned}
\end{equation}
Therefore,
\[
|\mathcal{U}(z',x_d)| \leq N x_d^{\frac{1-\alpha}{2}} C_2.
\]
Next, by Lemmas \ref{lem:caccioppoli} and \ref{lem:caccioimproved},
\begin{equation} \label{eq:6}
\begin{aligned}
&|u(z',x_d)| \leq |u(z',1/2)|+\int_{x_d}^{1/2}|D_d u(z',s)| ds \\
&\leq N C_1 + N(1+x_d^{\frac{1-\alpha}{2}-\delta})\left(\int_{0}^{1/2}\Vert D_d u(\cdot,s) \Vert^2_{W_2^{k/2,k}(Q_{1/2}')} s^{ \alpha} ds\right)^{1/2} \leq N(1+x_d^{\frac{1-\alpha}{2}-\delta})C_1.
\end{aligned}
\end{equation}
In the above, estimate $|u(z',1/2)| \leq N C_1$ follows from the interior $L_{\infty}$ estimate in Lemma \ref{lem:intbddnssest}, $\delta \in (0,1)$ is any constant when $\alpha \in \{1,3,5,\dots\}$, and $\delta=0$ otherwise. We introduce small $\delta>0$ to avoid the log-correction appearing in integrals $\int x_d^{-1} dx_d$. Applying estimate (\ref{eq:6}) to $u_t,D_{x'}u,D^2_{x'}u, \lambda u$, and $\lambda^{1/2}u$, and using Lemma \ref{lem:caccioimproved} give that on $Q_{1/2}^+$ we have
\begin{equation} \label{eq:esttangent1}
\begin{gathered}
|u_t|+|D_{x'}u|+|D^2_{x'}u|+\lambda|u|+\lambda^{1/2}|u| \leq N(1+x_d^{\frac{1-\alpha}{2}-\delta})C_2.
\end{gathered}
\end{equation}
Hence,
\begin{equation} \label{eq:estDd1}
|D_d u| \leq N|\mathcal{U}| + N |D_{x'}u| \leq N(1+x_d^{\frac{1-\alpha}{2}-\delta})C_2.
\end{equation}
From (\ref{eq:estDd1}) applied to $D_{x'} u$ and Lemma \ref{lem:caccioimproved}, we get
\begin{equation} \label{eq:estDdDtang1}
|D_d D_{x'}u| \leq N(1+x_d^{\frac{1-\alpha}{2}-\delta})C_2.
\end{equation}
Using (\ref{eq:esttangent1}) and (\ref{eq:estDdDtang1}) in the first line of (\ref{eq:estU1}), we get
\begin{equation*}
\begin{gathered}
|x_d^{\alpha}\mathcal{U}| \leq NC_2\int_0^{x_d}s^{\alpha}(1+s^{\frac{1-\alpha}{2}-\delta}) ds \leq N(x_d^{1+\alpha}+x_d^{\frac{3+\alpha}{2}-\delta})C_2,
\end{gathered}
\end{equation*}
or equivalently,
\[
|\mathcal{U}/x_d| \leq N(1+x_d^{\frac{1-\alpha}{2}-\delta})C_2.
\]

If $\alpha \in (-1,1)$, we are done with the case $p=2$. If $\alpha \geq 1$, we iterate the estimates above to improve the exponent $\frac{1-\alpha}{2}-\delta$ as follows. Plugging (\ref{eq:estDd1}) into the first line of (\ref{eq:6}), we get
\begin{equation*}
|u| \leq N(1+x_d^{\frac{3-\alpha}{2}-\delta})C_1.
\end{equation*}
Applying this estimate to $u_t,D_{x'}u,D^2_{x'}u, \lambda u$, and $\lambda^{1/2}u$, and using Lemma \ref{lem:caccioimproved}, we have
\[
|u_t|+|D_{x'}u|+|D^2_{x'}u|+\lambda|u|+\lambda^{1/2}|u| \leq N(1+x_d^{\frac{3-\alpha}{2}-\delta})C_2.
\]
Then,
\[
|D_d u| + |D_dD_{x'}u| \leq N(1+x_d^{\frac{3-\alpha}{2}-\delta})C_2.
\]
Using the above estimates in the first line of (\ref{eq:estU1}) gives
\[
|\mathcal{U}/x_d| \leq N(1+x_d^{\frac{3-\alpha}{2}-\delta})C_2.
\]
If $\alpha \in [1,3)$ then we are done. If $\alpha \geq 3$, by again iterating these estimates, in finitely many steps we arrive to
\begin{equation*}
|u| \leq NC_1, \quad
|Du| +\lambda^{1/2}|u| \leq NC_2, \quad
|\mathcal{U}/x_d| \leq N C_2.
\end{equation*}

For general $\tau \in [1/2,1)$, we cover $Q_{\tau r}^+$ with boundary cylinders of the form $Q_{(1-\tau)r}^+(z_0',0)$ with $z_0' \in \overline{Q_{\tau r}'}$, and interior cylinders $Q_{(1-\tau)r/4}(z_0)$ with $z_0 \in \overline{Q_{\tau r}^+}$ having $x_{0d} \geq (1-\tau)r/2$. On each boundary cylinder we apply the estimates proven above for $\tau=1/2$, and on each interior cylinder we use the estimates from Lemma \ref{lem:intbddnssest}. We then use inequality (\ref{eq:mudoubling1d}) for the measure $\mu$ to bound averages over the boundary and interior cylinders by the averages over $Q_r^+$. This completes the proof in the case $p=2$. Then by a standard iteration argument, see \cite[p. 184-186]{GiaquintaMartinazzi}, estimates (\ref{eq:bdrybddnessestsimplehom}) hold for general $p \in (1,2]$, assuming that initially $u \in \mathcal{H}^1_2(Q_r^+, \mu)$.

Assume now that $u \in \mathcal{H}^1_p(Q_r^+, \mu)$, where $p \in (1,2)$, and let $\tau \in [1/2,1)$ be fixed. Let $u^{(\varepsilon)}$ be the mollification of $u$ in the $t,x'$ variables. Let $r_1 = \frac{1+\tau}{2} r$ and $r_2 = \frac{3+\tau}{4} r$. Note that the argument above can essentially be repeated by taking
\[
C_1 = C_2 = \left(\int_0^{r_2} \Bigl[ \Vert Du^{(\varepsilon)}(\cdot, s) \Vert_{W^{k/2,k}_p(Q_{r_2}')}^p + \Vert u^{(\varepsilon)}(\cdot,s) \Vert_{W^{k/2,k}_p(Q_{r_2}')}^p \Bigr] s^{\alpha}ds \right)^{1/p} < \infty,
\]
where $k=k(p,d)$ is a fixed large enough even number. Indeed, by following the derivation of estimates above, the following sequence of estimates can be obtained for $z \in Q_{r_1}^+$. Let $\delta \in (0,1)$ be any constant if $\frac{1+\alpha}{p} \in \{1,2,3,\dots\}$ and $\delta=0$ otherwise. Then
\begin{gather*}
|\mathcal{U}^{(\varepsilon)}| \lesssim x_d^{1-\frac{1+\alpha}{p}}, \\
|u^{(\varepsilon)}| + |Du^{(\varepsilon)}|+|\mathcal{U}^{(\varepsilon)}/x_d| \lesssim 1 + x_d^{1-\frac{1+\alpha}{p}-\delta}, \\
|u^{(\varepsilon)}| + |Du^{(\varepsilon)}|+|\mathcal{U}^{(\varepsilon)}/x_d| \lesssim 1 + x_d^{2-\frac{1+\alpha}{p}-\delta},
\end{gather*}
and so on. In finitely many steps, we get $|u^{(\varepsilon)}|+|Du^{(\varepsilon)}|+ |\mathcal{U}^{(\varepsilon)}/x_d| \in L_{\infty}(Q_{r_1}^+)$.  In particular, $u^{(\varepsilon)} \in \mathcal{H}^1_2(Q_{r_1}^+, \mu)$ for all small enough $\varepsilon$. Then, (\ref{eq:bdrybddnessestsimplehom}) holds with $u^{(\varepsilon)}$ in place of $u$ in the left hand sides. Passing to the limit along a subsequence as $\varepsilon \to 0^+$, we obtain (\ref{eq:bdrybddnessestsimplehom}) for $u$.
\end{proof}

\begin{remark}
It should be noted that in the scalar case $n=1$ estimates (\ref{eq:bdrybddnessestsimplehom}) were proved in \cite[pp. 1481--1487]{DongPhanIndiana2023} with a somewhat different approach. There the authors first used Moser's iteration to show the bound $\Vert u \Vert_{L_{\infty}(Q_{1/2}^+)} \leq N\Vert u \Vert_{L_2(Q_1^+,\mu)}$. They then applied this estimate in an iteration scheme similar to the one we presented above to show the boundedness of $D_d u$ and $\mathcal{U}/x_d$. Due to the use of Moser's iteration, their argument cannot be adapted to the case of systems with $n \geq 2$.
\end{remark}

For $r \in (0,1]$ and $Q \subset \R^{d+1}$, we define the $(r/2, r)-$H\"older semi-norm of a function $f:Q\to X$, where $X$ is a Banach space, by
\[
[f]_{C^{r/2,r}(Q)} = \sup_{(t,x),(s,y) \in Q, (t,x) \neq (s,y)} \frac{|f(t,x)-f(s,y)|_X}{|t-s|^{r/2}+|x-y|^r}.
\]
From Lemma \ref{lem:bdrbddnssest} we obtain the following boundary Lipschitz estimates of solutions.
\begin{lem} \label{lem:bdrSch}
Under the assumptions of Lemma {\rm \ref{lem:bdrbddnssest}}, it holds that
\begin{equation*}
\begin{gathered}
[D_{x'}u]_{C^{1/2,1}(Q_{r/2}^+)}+[\mathcal{U}]_{C^{1/2,1}(Q_{r/2}^+)} + \sqrt{\lambda}[u]_{C^{1/2,1}(Q_{r/2}^+)}
\leq N r^{-1} \left( \fint_{Q_r^+} (|Du|^p+\lambda^{p/2} |u|^p) d\mu \right)^{1/p},
\end{gathered}
\end{equation*}
where $N=N(d,\alpha,p,\kappa) > 0$.
\end{lem}

\begin{proof}
By scaling, it is sufficient to consider $r=1$. Then
\begin{align*}
[D_{x'}u]_{C^{1/2,1}(Q_{1/2}^+)} \leq N(\Vert DD_{x'} u\Vert_{L_{\infty}(Q_{1/2}^+)}+\Vert \partial_t D_{x'} u\Vert_{L_{\infty}(Q_{1/2}^+)}) 
\leq N \left( \fint_{Q_1^+} (|Du|^p+\lambda^{p/2} |u|^p) d\mu \right)^{1/p},
\end{align*}
where the last inequality follows by Lemmas \ref{lem:bdrbddnssest} and \ref{lem:caccioimproved} applied to $D_{x'} u$ and $\partial_t D_{x'} u$. The terms $[\mathcal{U}]_{C^{1/2,1}(Q_{1/2}^+)}$ and $\sqrt{\lambda}[u]_{C^{1/2,1}(Q_{1/2}^+)}$ are estimated similarly. When estimating $[\mathcal{U}]_{C^{1/2,1}(Q_{1/2}^+)}$, we use the estimate
$
|D_d \mathcal{U}| \leq N(|u_t|+|DD_{x'}u|+\lambda|u|+|\mathcal{U}/x_d|)
$, which follows from the equation.
\end{proof}

Similarly, we have interior Lipschitz estimates of solutions.
\begin{lem} \label{lem:intSch}
Under the assumptions of Lemma {\rm \ref{lem:intbddnssest}}, we have
\begin{equation*}
\begin{gathered}
[D_{x'}u]_{C^{1/2,1}(Q_{r/2}(z_0))}+[\mathcal{U}]_{C^{1/2,1}(Q_{r/2}(z_0))} + \sqrt{\lambda}[u]_{C^{1/2,1}(Q_{r/2}(z_0))} \\
\leq N r^{-1} \left( \fint_{Q_r(z_0)} (|Du|^p+\lambda^{p/2} |u|^p) d\mu \right)^{1/p},
\end{gathered}
\end{equation*}
where $N=N(d,\alpha,p,\kappa) > 0$.
\end{lem}

From Lemmas \ref{lem:bdrSch} and \ref{lem:intSch}, we get the following mean oscillation estimates.
\begin{lem} \label{lem:homoscest}
Let $p \in (1,2]$, $z_0 \in \overline{\R^{d+1}_+}, r \in (0,\infty), \nu \in [32,\infty)$, and $u \in \mathcal{H}^1_p(Q_{\nu r}^+(z_0), \mu; \R^n)$ be a weak solution to \eqref{eq:simplehom}-\eqref{eq:conormalsimplehom} in $Q_{\nu r}^+(z_0)$. Then
\begin{align*}
&\fint_{Q_r^+(z_0)}|D_{x'}u-(D_{x'}u)_{Q_r^+(z_0), \mu}| d\mu+\fint_{Q_r^+(z_0)}|\mathcal{U}-(\mathcal{U})_{Q_r^+(z_0), \mu}| d\mu \\
&+\sqrt{\lambda}\fint_{Q_r^+(z_0)}|u-(u)_{Q_r^+(z_0), \mu}| d\mu
\leq N \nu^{-1} \left( \fint_{Q^+_{\nu r}(z_0)} (|Du|^p+\lambda^{p/2} |u|^p) d\mu \right)^{1/p},
\end{align*}
where $N=N(d,\alpha,p,\kappa) > 0$.
\end{lem}

\begin{proof}
Case I: $x_{0d} \geq \frac{\nu r}{8} = 2 \frac{\nu r}{16}$. Then
\begin{align*}
&\fint_{Q_r^+(z_0)}|D_{x'}u-(D_{x'}u)_{Q_r^+(z_0), \mu}| d\mu \leq Nr[D_{x'}u]_{C^{1/2,1}(Q_{\frac{\nu r}{32}}(z_0))} \leq [\text{by Lemma \ref{lem:intSch}}] \\
&\leq N \nu^{-1} \left( \fint_{Q_{\frac{\nu r}{16}}(z_0)} (|Du|^p+\lambda^{p/2} |u|^p) d\mu \right)^{1/p} \leq N \nu^{-1} \left( \fint_{Q^+_{\nu r}(z_0)} (|Du|^p+\lambda^{p/2} |u|^p) d\mu \right)^{1/p}.
\end{align*}
Case II: $x_{0d} < \frac{\nu r}{8}$. Let $\Hat{z}_0 = (t_0,x_0',0)$. Then
\begin{align*}
&\fint_{Q_r^+(z_0)}|D_{x'}u-(D_{x'}u)_{Q_r^+(z_0), \mu}| d\mu \leq N r[D_{x'}u]_{C^{1/2,1}(Q^+_{\frac{\nu r}{4}}(\Hat{z}_0))} \leq [\text{by Lemma \ref{lem:bdrSch}}] \\
&\leq N \nu^{-1} \left( \fint_{Q^+_{\frac{\nu r}{2}}(\Hat{z}_0)} (|Du|^p+\lambda^{p/2} |u|^p) d\mu \right)^{1/p} \leq N \nu^{-1} \left( \fint_{Q^+_{\nu r}(z_0)} (|Du|^p+\lambda^{p/2} |u|^p) d\mu \right)^{1/p}.
\end{align*}
For $\mathcal{U}$ and $\sqrt{\lambda} u$, we use the same argument as above.
\end{proof}

\subsection{\texorpdfstring{$L_p$}{} estimates for inhomogeneous systems}
Before we give the proof of Theorem {\rm \ref{thm:exisuniqsimpleunweighted}}, some more preparation is needed. We now prove the following result on decomposition of solutions.
\begin{proposition} \label{prop:decompsimple}
Let $z_0 \in \overline{\R^{d+1}_+}, r>0$, $\lambda\ge 0$, and $\nu \in [32,\infty)$. Suppose that $G = |F|+|f| \in L_2(Q_{\nu r}^+(z_0), \mu)$ and $u \in \mathcal{H}^1_2(Q_{\nu r}^+(z_0), \mu; \R^n)$ is a weak solution to \eqref{eq:simplenonhom}-\eqref{eq:conormalsimplenonhom} in $Q_{\nu r}^+(z_0)$. Then, we can write $u(t,x)=v(t,x)+w(t,x)$ in $Q_{\nu r}^+(z_0)$, where $v$ and $w$ are in $\mathcal{H}^1_2(Q_{\nu r}^+(z_0),\mu; \R^n)$ and satisfy
\begin{equation}
\label{eq:decompV}
\fint_{Q_{\nu r}^+(z_0)}|V|^2 d\mu \leq N \fint_{Q_{\nu r}^+(z_0)} |G|^2 d\mu,
\end{equation}
\begin{equation}
\label{eq:decompW}
\begin{aligned}
&\fint_{Q_r^+(z_0)}|D_{x'}w-(D_{x'}w)_{Q_r^+(z_0), \mu}| d\mu+\fint_{Q_r^+(z_0)}|\mathcal{W}-(\mathcal{W})_{Q_r^+(z_0), \mu}| d\mu \\
&\,\,+\sqrt{\lambda}\fint_{Q_r^+(z_0)}|w-(w)_{Q_r^+(z_0), \mu}| d\mu \leq N \nu^{-1}\left(\fint_{Q_{\nu r}^+(z_0)} |U|^2 d\mu\right)^{1/2} + N\nu^{-1}\left(\fint_{Q_{\nu r}^+(z_0)} |G|^2 d\mu\right)^{1/2},
\end{aligned}
\end{equation}
where $N=N(d,\alpha, \kappa) > 0$ and
\[
V = |Dv|+\lambda^{1/2}|v|, \quad U = |Du|+\lambda^{1/2}|u|, \quad \mathcal{W} = A^{dj}D_j w.
\]
\end{proposition}
\begin{proof}
We may assume that $\lambda>0$. The case when $\lambda=0$ follows by adding $\eta x_d^{\alpha}u$ with $\eta>0$ to both sides of the equation and then taking the limit as $\eta\to 0$. Let $v \in \mathcal{H}^1_2(\Omega_{t_0}, \mu; \R^n)$ be the weak solution to
\begin{equation*}
\begin{cases}
\mathcal{P}_0v = -D_i(x_d^{\alpha}F^i\chi_{Q_{\nu r}^+(z_0)})+ \sqrt{\lambda}x_d^{\alpha} f\chi_{Q_{\nu r}^+(z_0)} \quad \text{in $\Omega_{t_0}$}, \\
\lim_{x_d \to 0^+} x_d^{\alpha}(A^{dj}(x_d)D_j v - F^d \chi_{Q_{\nu r}^+(z_0)}) = 0.
\end{cases}
\end{equation*}
Then, (\ref{eq:decompV}) follows by Theorem \ref{thm:L2solv}.

Let $w=u-v$. Then $w \in \mathcal{H}^1_2(Q_{\nu r}^+(z_0), \mu; \R^n)$ is a weak solution to the homogeneous system
\[
\mathcal{P}_0w = 0 \quad \text{in $Q_{\nu r}^+(z_0)$}
\]
with the conormal boundary condition. Estimate (\ref{eq:decompW}) now follows from Lemma \ref{lem:homoscest} and the triangle inequality.
\end{proof}

In the following two lemmas, we prove that $\mathcal{H}^1_2$ solutions also belong to $\mathcal{H}^1_p$ provided that $F$ and $f$ are sufficiently regular.
\begin{lem} \label{lem:LptoLsmall}
Let $\lambda > 0$ and $F^1,\dots,F^d,f \in L_{\infty}(\R^{d+1}_+; \R^n)$ have bounded supports. Let $u \in \mathcal{H}^1_2(\R^{d+1}_+, \mu; \R^n)$ be the weak solution to \eqref{eq:simplenonhom}-\eqref{eq:conormalsimplenonhom} in $\R^{d+1}_+$. Then $u \in \mathcal{H}^1_p(\R^{d+1}_+, \mu)$ for any $p \in (1,2)$.
\end{lem}

\begin{proof}
Let $(R_k)_{k=0}^{\infty}$ be an increasing sequence of positive numbers with $R_k \to \infty$ as $k \to \infty$, which is to be specified later. Denote $\Hat{Q}_k = (-R_k^2,R_k^2) \times D_{R_k}^+$. Fix $\eta_k \in C^{\infty}(\R^{d+1})$ such that $\eta_k = 1$ on $\R^{d+1}_+ \setminus \Hat{Q}_{k+1}$, $\eta_k = 0$ on $\Hat{Q}_k$, $|D \eta_k| \leq \frac{N}{R_{k+1}-R_k}$, and $|\partial_t \eta_k| \leq \frac{N}{(R_{k+1}-R_k)^2}$. Take $R_0$ large enough so that $\operatorname{supp}F \cup \operatorname{supp}f \subset \overline{\Hat{Q}_0}$.

By H\"older's inequality
\[
\Vert Du \Vert_{L_p(\Hat{Q}_1,\mu)} + \Vert u \Vert_{L_p(\Hat{Q}_1,\mu)} < \infty.
\]
Hence it remains to prove that
\[
\Vert Du \Vert_{L_p(\R^{d+1}_+\setminus \Hat{Q}_0,\mu)} + \Vert u \Vert_{L_p(\R^{d+1}_+\setminus \Hat{Q}_0,\mu)} < \infty.
\]
Let $u_k = u \eta_k$. Then $u_k$ is a solution to
\[
\begin{cases}
\mathcal{P}_0 u_k = -D_i(x_d^{ \alpha} F_{k}^i) + \sqrt{\lambda}f_k \quad \text{in $\R^{d+1}_+$}, \\
\lim_{x_d \to 0^+} x_d^{\alpha}(A^{dj}D_j u_k - F_{k}^d) = 0,
\end{cases}
\]
where
\[
F_{k}^i = A^{ij} u D_j \eta_k,\quad 
f_k = \lambda^{-1/2}(-A^{ij}D_j u D_i \eta_k + a_0 u \partial_t\eta_k).
\]
From the $L_2$ estimate (Theorem \ref{thm:L2solv}), we have
\begin{gather*}
\Vert Du\Vert_{L_2(\Hat{Q}_{k+2}\setminus \Hat{Q}_{k+1},\mu)} + \lambda^{1/2} \Vert u\Vert_{L_2(\Hat{Q}_{k+2}\setminus \Hat{Q}_{k+1},\mu)} \leq \Vert D(u\eta_k)\Vert_{L_2(\R^{d+1}_+, \mu)} + \lambda^{1/2} \Vert u \eta_k\Vert_{L_2(\R^{d+1}_+, \mu)} \\
\leq N \Vert F_k\Vert_{L_2(\R^{d+1}_+, \mu)} + N \Vert f_k \Vert_{L_2(\R^{d+1}_+, \mu)} \leq N \lambda^{-1/2}(R_{k+1}-R_k)^{-1}\left(\Vert Du \Vert_{L_2(\Hat{Q}_{k+1}\setminus \Hat{Q}_k,\mu)} \right. \\
\left. + \lambda^{1/2}\Vert u \Vert_{L_2(\Hat{Q}_{k+1}\setminus \Hat{Q}_k,\mu)}+ (R_{k+1}-R_k)^{-1}\Vert u \Vert_{L_2(\Hat{Q}_{k+1}\setminus \Hat{Q}_k,\mu)}\right).
\end{gather*}
Take $R_{k+1}-R_k = C \geq \lambda^{-1/2}$, where $C$ is a constant to be specified later. Then
\begin{multline*}
\Vert Du\Vert_{L_2(\Hat{Q}_{k+2}\setminus \Hat{Q}_{k+1},\mu)} + \lambda^{1/2} \Vert u\Vert_{L_2(\Hat{Q}_{k+2}\setminus \Hat{Q}_{k+1},\mu)} 
\\ \leq N_0 \lambda^{-1/2}C^{-1}(\Vert Du \Vert_{L_2(\Hat{Q}_{k+1}\setminus \Hat{Q}_k,\mu)} + \lambda^{1/2}\Vert u \Vert_{L_2(\Hat{Q}_{k+1}\setminus \Hat{Q}_k,\mu)})
\end{multline*}
for some $N_0>0$ which we now fix through the end of the proof. By iterating the above estimate,
\begin{multline*}
\Vert Du \Vert_{L_2(\Hat{Q}_{k+1}\setminus \Hat{Q}_k,\mu)} + \lambda^{1/2}\Vert u \Vert_{L_2(\Hat{Q}_{k+1}\setminus \Hat{Q}_k,\mu)} \\
\leq (N_0 \lambda^{-1/2} C^{-1})^k (\Vert Du \Vert_{L_2(\Hat{Q}_{1}\setminus \Hat{Q}_0,\mu)} + \lambda^{1/2}\Vert u \Vert_{L_2(\Hat{Q}_{1}\setminus \Hat{Q}_0,\mu)}).
\end{multline*}
By H\"older's inequality, we get
\begin{gather*}
\Vert Du \Vert_{L_p(\Hat{Q}_{k+1}\setminus \Hat{Q}_k,\mu)} + \lambda^{1/2}\Vert u \Vert_{L_p(\Hat{Q}_{k+1}\setminus \Hat{Q}_k,\mu)} \leq (\mu(\Hat{Q}_{k+1}))^{\frac{1}{p}-\frac{1}{2}} \left(\Vert Du \Vert_{L_2(\Hat{Q}_{k+1}\setminus \Hat{Q}_k,\mu)} \right.\\
\left. + \lambda^{1/2}\Vert u \Vert_{L_2(\Hat{Q}_{k+1}\setminus \Hat{Q}_k,\mu)}\right) \leq \Bigl[\mu(\Hat{Q}_{k+1}) = N_{d,\alpha} R_{k+1}^{d+2+\alpha} = N_{d,\alpha} (R_0+(k+1)C)^{d+2+\alpha}\Bigr] \\
\leq N (R_0+(k+1)C)^{(d+2+\alpha)(\frac{1}{p}-\frac{1}{2})} (N_0 \lambda^{-1/2} C^{-1})^k (\Vert Du \Vert_{L_2(\Hat{Q}_{1}\setminus \Hat{Q}_0,\mu)} + \lambda^{1/2}\Vert u \Vert_{L_2(\Hat{Q}_{1}\setminus \Hat{Q}_0,\mu)}).
\end{gather*}
Taking $C = 2.1 N_0 \lambda^{-1/2}$, we get
\[
\Vert Du \Vert_{L_p(\Hat{Q}_{k+1}\setminus \Hat{Q}_k,\mu)} + \lambda^{1/2}\Vert u \Vert_{L_p(\Hat{Q}_{k+1}\setminus \Hat{Q}_k,\mu)} \leq N 2^{-k} (\Vert Du \Vert_{L_2(\Hat{Q}_{1}\setminus \Hat{Q}_0,\mu)} + \lambda^{1/2}\Vert u \Vert_{L_2(\Hat{Q}_{1}\setminus \Hat{Q}_0,\mu)}),
\]
where $N$ can now depend on $\lambda$. Then,
\begin{gather*}
\Vert Du \Vert_{L_p(\R^{d+1}_+\setminus \Hat{Q}_0,\mu)}^p + \lambda^{p/2}\Vert u \Vert_{L_p(\R^{d+1}_+\setminus \Hat{Q}_0,\mu)}^p = \sum_{k=0}^\infty (\Vert Du \Vert_{L_p(\Hat{Q}_{k+1}\setminus \Hat{Q}_k,\mu)}^p + \lambda^{p/2}\Vert u \Vert_{L_p(\Hat{Q}_{k+1}\setminus \Hat{Q}_k,\mu)}^p) \\
\leq N (\Vert Du \Vert_{L_2(\Hat{Q}_{1}\setminus \Hat{Q}_0,\mu)} + \lambda^{1/2}\Vert u \Vert_{L_2(\Hat{Q}_{1}\setminus \Hat{Q}_0,\mu)})^p \sum_{k=0}^\infty 2^{-kp} < \infty.
\end{gather*}
The lemma is proved.
\end{proof}

\begin{lem} \label{lem:L2toLinfty}
Let $\lambda > 0$ and $F^1,\dots,F^d,f \in C_0^{\infty}(\R^{d+1}_+;\R^n)$. Let $u \in \mathcal{H}^1_2(\R^{d+1}_+, \mu; \R^n)$ be the weak solution to \eqref{eq:simplenonhom}-\eqref{eq:conormalsimplenonhom} in $\R^{d+1}_+$. Then $u, Du \in L_{p}(\R^{d+1}_+, \mu)$ for any $p \in (2, \infty]$.
\end{lem}
\begin{proof} By interpolation, it is sufficient to prove the claim for $p = \infty$. With $F^i,f \in C_0^{\infty}(\R^{d+1}_+; \R^n)$, we can essentially repeat the proofs of Lemmas \ref{lem:intbddnssest} and \ref{lem:bdrbddnssest}, and therefore,
\[
h(z_0) := \Vert Du \Vert_{L_{\infty}(Q_1^+(z_0))} + \Vert u \Vert_{L_{\infty}(Q_1^+(z_0))} < \infty
\]
for any $z_0 \in \overline{\R^{d+1}_+}$. Moreover, $h(z_0) \to 0$ as $z_0 \to \infty$ by Lemmas \ref{lem:intbddnssest} and \ref{lem:bdrbddnssest}. Hence
\[
\Vert Du \Vert_{L_{\infty}(\R^{d+1}_+)} + \Vert u \Vert_{L_{\infty}(\R^{d+1}_+)} < \infty.
\]
The lemma is proved.
\end{proof}

The Hardy-Littlewood and Fefferman-Stein maximal functions of an, in general, a vector-valued function $f$ on $\Omega_T$ with respect to the measure $\mu$ are respectively defined by
\begin{gather*}
\mathcal{M} f (z) = \sup_{z \in Q_r^+(z_0)} \fint_{Q_r^+(z_0)} |f(y)| d\mu(y), \\
f^{\sharp}(z) = \sup_{z \in Q_r^+(z_0)} \fint_{Q_r^+(z_0)} |f(y)-(f)_{Q_r^+(z_0), \mu}| d\mu(y),
\end{gather*}
where the $\sup$ are taken over all cylinders $Q_r^+(z_0)$ containing $z$ with $z_0 \in \overline{\Omega_T}$ and $r>0$. The following versions of weighted Hardy-Littlewood and Fefferman-Stein inequalities are special cases of Corollaries 2.6 and 2.7 in \cite{DongKimTAMS2018}.

\begin{thm} \label{thm:HLFS}
Let $\alpha \in (-1,\infty), p,q \in (1,\infty), K \in [1,\infty)$, and $w=w_0(t)w_1(x)$, where $w_0 \in A_q(\R), w_1 \in A_p(\R^d_+,\mu)$ with $[w_0]_{A_q(\R)} \leq K, [w_1]_{A_p(\R^d_+,\mu)} \leq K$. Then, for any $f \in L_{q,p}(\Omega_T, wd\mu)$ holds
\begin{gather*}
\Vert \mathcal{M} f \Vert_{L_{q,p}(\Omega_T, wd\mu)} \leq N \Vert f \Vert_{L_{q,p}(\Omega_T, wd\mu)}, \\
\Vert f \Vert_{L_{q,p}(\Omega_T, wd\mu)} \leq N \Vert f^{\sharp} \Vert_{L_{q,p}(\Omega_T, wd\mu)},
\end{gather*}
where $N=N(d,\alpha,p,q,K)>0$.
\end{thm}

Recall that we use the notation $U = |Du| + \lambda^{1/2}|u|$.
\begin{proof}[Proof of Theorem {\rm \ref{thm:exisuniqsimpleunweighted}}]
When $p=2$, the theorem follows by Theorem \ref{thm:L2solv}.

We now prove the a-priori estimate (\ref{eq:apriorisimpleunweighted}) in the case $p \in (2,\infty)$. Let $\nu \in [32,\infty)$ be a constant which will be specified later. Let $z_0 \in \overline{\Omega_T}$ and $r>0$. By Proposition \ref{prop:decompsimple} and the triangle inequality, it holds that
\begin{equation*}
\begin{gathered}
\fint_{Q_r^+(z_0)}|D_{x'}u-(D_{x'}u)_{Q_r^+(z_0)}| d\mu+\fint_{Q_r^+(z_0)}|\mathcal{U}-(\mathcal{U})_{Q_r^+(z_0)}| d\mu \\
+\sqrt{\lambda}\fint_{Q_r^+(z_0)}|u-(u)_{Q_r^+(z_0)}| d\mu \leq N \nu^{-1}\left(\fint_{Q_{\nu r}^+(z_0)} |U|^2 d\mu\right)^{1/2} \\
+ N(\nu^{-1}+\nu^{\frac{d+2+\alpha_+}{2}})\left(\fint_{Q_{\nu r}^+(z_0)} |F|^2+|f|^2 d\mu\right)^{1/2},
\end{gathered}
\end{equation*}
where we used $\mu(Q_{\nu r}^+(z_0)) \leq N_{\alpha} \nu^{d+2+\alpha_+}\mu(Q_r^+(z_0))$. It now follows that on $\Omega_T$, we have
\begin{equation*}
(D_{x'} u)^{\sharp} + (\mathcal{U})^{\sharp} + (\sqrt{\lambda} u)^{\sharp} \leq N \nu^{-1} \mathcal{M}^{1/2}(|U|^2) + N\nu^{\frac{d+2+\alpha_+}{2}}(\mathcal{M}^{1/2}(|F|^2)+\mathcal{M}^{1/2}(|f|^2)).
\end{equation*}
By taking the $L_p(\Omega_T, \mu)$ norms in the above inequality, applying the weighted Hardy-Littlewood and Fefferman-Stein inequalities (Theorem \ref{thm:HLFS} with $q=p$ and $w=1$), and using $|D_d u| \leq N(|D_{x'} u| + |\mathcal{U}|)$, we have
\begin{gather*}
\Vert Du \Vert_{L_p(\Omega_T, \mu)} + \sqrt{\lambda} \Vert u\Vert_{L_p(\Omega_T, \mu)} \leq N\nu^{-1} (\Vert Du \Vert_{L_p(\Omega_T, \mu)} + \sqrt{\lambda} \Vert u\Vert_{L_p(\Omega_T, \mu)}) \\
+ N\nu^{\frac{d+2+\alpha_+}{2}} (\Vert F \Vert_{L_p(\Omega_T, \mu)} + \Vert f\Vert_{L_p(\Omega_T, \mu)}).
\end{gather*}
Taking $\nu$ large enough ($\nu = \max(32, 2N)$), we arrive to (\ref{eq:apriorisimpleunweighted}).

To prove the existence of solutions with $p \in (2,\infty)$, assume without loss that $T=\infty$, that is $\Omega_T = \R^{d+1}_+$, and take sequences $F_k, f_k \in C_0^{\infty}(\R^{d+1}_+)$ such that $F_k \to F$ and $f_k \to f$ in $L_p(\R^{d+1}_+,\mu)$. Then by the $L_2$ solvability, there exist solutions $u_k \in \mathcal{H}^1_2(\R^{d+1}_+,\mu; \R^n)$ to \eqref{eq:simplenonhom}-\eqref{eq:conormalsimplenonhom} with $F_k$ and $f_k$ in place of $F$ and $f$. By Lemma \ref{lem:L2toLinfty}, $u_k \in \mathcal{H}^1_p(\R^{d+1}_+,\mu)$. By the a-priori estimate (\ref{eq:apriorisimpleunweighted}), the sequence $u_k$ is Cauchy in $\mathcal{H}^1_p(\R^{d+1}_+,\mu)$, and therefore converges to the solution $u \in \mathcal{H}^1_p(\R^{d+1}_+,\mu; \R^n)$ with $F,f$.

For $p \in (1,2)$, we use a duality argument. Let  $G \in L_{p'}(\Omega_T, \mu; \R^n)^d$ and $g \in L_{p'}(\Omega_T, \mu; \R^n)$ be arbitrary. If $T<\infty$, we also extend $G$ and $g$ to be zero for $t \geq T$.  As $p' \in (2,\infty)$, by reversing the time variable, there exists a unique weak solution $v \in \mathcal{H}^1_{p'}(\R^{d+1}_+, \mu; \R^n)$ to the adjoint equation
\[
\begin{cases}
-x_d^{\alpha}a_0(x_d)v_t -D_i(x_d^{\alpha}(A^{ji})^T D_j v) + \lambda x_d^{\alpha} C_0^T v= -D_i(x_d^{\alpha} G^i) + \sqrt{\lambda}x_d^{\alpha}g \quad \text{in $\R^{d+1}_+$}, \\
\lim_{x_d \to 0^+} x_d^{\alpha}((A^{jd})^T D_j v - G^d) = 0.
\end{cases}
\]
Moreover, $v(t,\cdot)=0$ for $t \geq T$ if $T < \infty$, and
\[
\Vert Dv \Vert_{L_{p'}(\Omega_T, \mu)} + \sqrt{\lambda}\Vert v \Vert_{L_{p'}(\Omega_T, \mu)} \leq N \Vert G \Vert_{L_{p'}(\Omega_T, \mu)} + N \Vert g \Vert_{L_{p'}(\Omega_T, \mu)}.
\]
By taking $\varphi = v$ in (\ref{eq:soldefvar}) and using the equations satisfied by $u$ and $v$, we get
\begin{gather*}
\int_{\Omega_T} (G^i \cdot D_i u + g\cdot \sqrt{\lambda}u)d\mu = \int_{\Omega_T} (F^i \cdot D_i v + f\cdot \sqrt{\lambda}v)d\mu.
\end{gather*}
Here one can use Steklov averages to justify taking $u$ and $v$ as test functions in the weak formulation of solutions. We have
\begin{gather*}
\left|\int_{\Omega_T} (G^i \cdot D_i u + g\cdot \sqrt{\lambda}u) d\mu \right|
\leq N (\Vert F\Vert_{L_{p}(\Omega_T,\mu)} + \Vert f\Vert_{L_{p}(\Omega_T,\mu)})(\Vert G \Vert_{L_{p'}(\Omega_T, \mu)} + \Vert g \Vert_{L_{p'}(\Omega_T, \mu)}).
\end{gather*}
Since $G$ and $g$ are arbitrary, estimate (\ref{eq:apriorisimpleunweighted}) follows.

To prove the existence of solutions with $p \in (1,2)$, again assume without loss that $T=\infty$, and take sequences $F_k, f_k \in L_{\infty}(\R^{d+1}_+)$ with compact supports such that $F_k \to F$ and $f_k \to f$ in $L_p(\R^{d+1}_+,\mu)$. Then by the $L_2$ solvability, there exist solutions $u_k \in \mathcal{H}^1_2(\R^{d+1}_+,\mu; \R^n)$ to \eqref{eq:simplenonhom}-\eqref{eq:conormalsimplenonhom} with $F_k$ and $f_k$ in place of $F$ and $f$. By Lemma \ref{lem:LptoLsmall}, $u_k \in \mathcal{H}^1_p(\R^{d+1}_+,\mu)$. By the a-priori estimate (\ref{eq:apriorisimpleunweighted}), the sequence $u_k$ is Cauchy in $\mathcal{H}^1_p(\R^{d+1}_+,\mu)$, and therefore converges to the solution $u \in \mathcal{H}^1_p(\R^{d+1}_+,\mu; \R^n)$ to \eqref{eq:simplenonhom}-\eqref{eq:conormalsimplenonhom}  with $F,f$.
\end{proof}

With Theorem \ref{thm:exisuniqsimpleunweighted} now available, we obtain an improved version of Proposition \ref{prop:decompsimple}, which we will need later.
\begin{proposition} \label{prop:decompsimplegenp}
Let $p \in (1,\infty), z_0 \in \overline{\R^{d+1}_+}, r>0$, and $\nu \in [32,\infty)$. Suppose that $G = |F|+|f| \in L_p(Q_{\nu r}^+(z_0), \mu)$ and $u \in \mathcal{H}^1_p(Q_{\nu r}^+(z_0), \mu; \R^n)$ is a weak solution to \eqref{eq:simplenonhom}-\eqref{eq:conormalsimplenonhom} in $Q_{\nu r}^+(z_0)$. Then, we can write $u(t,x)=v(t,x)+w(t,x)$ in $Q_{\nu r}^+(z_0)$, where $v$ and $w$ are in $\mathcal{H}^1_p(Q_{\nu r}^+(z_0),\mu; \R^n)$ and satisfy
\begin{equation*} 
\begin{aligned}
&\fint_{Q_{\nu r}^+(z_0)}|V|^p d\mu \leq N \fint_{Q_{\nu r}^+(z_0)} |G|^p d\mu, \\
&\fint_{Q_r^+(z_0)}|D_{x'}w-(D_{x'}w)_{Q_r^+(z_0)}| d\mu+\fint_{Q_r^+(z_0)}|\mathcal{W}-(\mathcal{W})_{Q_r^+(z_0)}| d\mu \\
&\quad +\sqrt{\lambda}\fint_{Q_r^+(z_0)}|w-(w)_{Q_r^+(z_0)}| d\mu \leq N \nu^{-1}\left(\fint_{Q_{\nu r}^+(z_0)} |U|^p d\mu\right)^{1/p} + N\nu^{-1}\left(\fint_{Q_{\nu r}^+(z_0)} |G|^p d\mu\right)^{1/p},
\end{aligned}
\end{equation*}
where $N=N(d,\alpha, p, \kappa)>0$ and $V, U, \mathcal{W}$ are as in Proposition {\rm \ref{prop:decompsimple}}.
\end{proposition}
The proof follows word-by-word the proof of Proposition \ref{prop:decompsimple} by changing $2$ to $p$ and using Theorem \ref{thm:exisuniqsimpleunweighted} instead of Theorem \ref{thm:L2solv}.

\section{Equations with partially BMO coefficients in weighted spaces} \label{sec:5}

In this section, we give the proof of Theorem \ref{thm:main1} (i). For now, we consider systems
\begin{equation} \label{eq:eqpar2}
\begin{cases}
x_d^{\alpha} a_0(x_d) u_t- D_i(x_d^{\alpha}A^{ij}(t,x)D_j u) + \lambda x_d^{\alpha} C_0(t,x) u= -D_i(x_d^{\alpha} F^i)+\sqrt{\lambda}x_d^{\alpha} f, \\
\lim_{x_d \to 0^+} x_d^{\alpha}(A^{dj}(t,x)D_j u - F^d) = 0.
\end{cases}
\end{equation}

\begin{proposition} \label{prop:prop3}
Let $z_0 \in \overline{\R^{d+1}_+}, \gamma_0 > 0, \nu \geq 32, R_0 \in (0,\infty]$, $r \in (0,\frac{R_0}{\nu}]$ $(r \in (0,\infty)$ if $R_0 = \infty)$, and $q_0, \sigma \in (1,\infty)$. Let $F^1,\dots,F^d,f \in L_{q_0}(Q_{\nu r}^+(z_0), \mu; \R^n)$ and $u \in \mathcal{H}^1_{q_0\sigma}(Q_{\nu r}^+(z_0), \mu; \R^n)$ be a solution to \eqref{eq:eqpar2} in $Q_{\nu r}^+(z_0)$. Then under Assumption {\rm \ref{assump} ($\gamma_0,R_0$)} there exists $U^Q:Q_{\nu r}^+(z_0) \to \R_{\geq 0}$ such that
\[
N^{-1}_{\kappa} U \leq U^Q \leq N_{\kappa}U \quad \text{on $Q_{\nu r}^+(z_0)$}, 
\]
and
\begin{equation} \label{eq:oscestU^Qgen}
\begin{gathered}
\fint_{Q_r^+(z_0)}|U^Q-(U^Q)_{Q_r^+(z_0), \mu}| d\mu \leq N \nu^{\frac{d+2+\alpha_+}{q_0}}\left(\fint_{Q_{\nu r}^+(z_0)}(|F|^{q_0}+|f|^{q_0})d\mu\right)^{1/{q_0}} \\
+ N(\nu^{-1}+\nu^{\frac{d+2+\alpha_+}{{q_0}}}\gamma_0^{\frac{1}{{q_0}\sigma'}})\left(\fint_{Q_{\nu r}^+(z_0)} |U|^{{q_0}\sigma} d\mu\right)^{\frac{1}{{q_0}\sigma}},
\end{gathered}
\end{equation}
where $N=N(d,\alpha,\kappa,q_0,\sigma) > 0$.
\end{proposition}
\begin{proof}
Let 
\[
\overline{A}^{ij}(x_d) = [A^{ij}]_{\nu r, z_0}(x_d), \quad \overline{C}_{0}(x_d) = [C_{0}]_{\nu r, z_0}(x_d)
\]
be from Assumption \ref{assump} ($\gamma_0,R_0$). Then $u$ satisfies in $Q_{\nu r}^+(z_0)$
\begin{equation*}
\begin{aligned}
&x_d^{\alpha}a_0(x_d) u_t - D_i(x_d^{\alpha}\overline{A}^{ij}(x_d)D_j u) + \lambda x_d^{\alpha} \overline{C}_0(x_d) u 
\\
&= -D_i(x_d^{\alpha}(F^i-(A^{ij}(z)-\overline{A}^{ij}(x_d))D_j u)) + \sqrt{\lambda}x_d^{\alpha} (f-\sqrt{\lambda}(C_0(z)-\overline{C}_0(x_d))u)
\end{aligned}
\end{equation*}
with the conormal boundary condition. By Proposition \ref{prop:decompsimplegenp}, we can write $u = v + w$ in $Q_{\nu r}^+(z_0)$, where
\begin{equation*}
\begin{aligned}
&\nu^{-\frac{d+2+\alpha_+}{q_0}}\left(\fint_{Q_{r}^+(z_0)}|V|^{q_0} d\mu\right)^{1/q_0} \leq N\left(\fint_{Q_{\nu r}^+(z_0)}|V|^{q_0} d\mu\right)^{1/q_0} \\
&\leq N \Bigl(\fint_{Q_{\nu r}^+(z_0)}(|F|^{q_0}+|f|^{q_0} + |A(z)-\overline{A}(x_d)|^{q_0} |Du|^{q_0} + |C_0(z)-\overline{C}_0(x_d)|^{q_0} |\sqrt{\lambda} u|^{q_0}) d\mu\Bigr)^{1/{q_0}} \\
&\leq [\text{by H\"older's inequality}] \leq N \left(\fint_{Q_{\nu r}^+(z_0)}(|F|^{q_0}+|f|^{q_0})d\mu\right)^{1/{q_0}} + N \gamma_0^{\frac{1}{q_0\sigma'}}\left(\fint_{Q_{\nu r}^+(z_0)}|U|^{q_0\sigma} d\mu\right)^{\frac{1}{q_0\sigma}}.
\end{aligned}
\end{equation*}
Here $V$ and $\mathcal{W}$ are as in Proposition \ref{prop:decompsimple}. Again by Proposition \ref{prop:decompsimplegenp}, we have
\begin{equation*}
\begin{aligned}
&\fint_{Q_r^+(z_0)}|D_{x'}w-(D_{x'}w)_{Q_r^+(z_0),\mu}| d\mu+\fint_{Q_r^+(z_0)}|\mathcal{W}-(\mathcal{W})_{Q_r^+(z_0),\mu}| d\mu+\sqrt{\lambda}\fint_{Q_r^+(z_0)}|w-(w)_{Q_r^+(z_0),\mu}| d\mu \\
&\leq N\nu^{-1}\left(\fint_{Q_{\nu r}^+(z_0)}(|F|^{q_0}+|f|^{q_0}) d\mu\right)^{1/{q_0}} 
+ N\nu^{-1}(1+\gamma_0^{\frac{1}{q_0\sigma'}})\left(\fint_{Q_{\nu r}^+(z_0)} |U|^{q_0\sigma} d\mu\right)^{\frac{1}{q_0\sigma}}.
\end{aligned}
\end{equation*}
Let $U^Q = |D_{x'}u|+\lambda^{1/2}|u|+|\overline{A}^{dj} D_j u|$. Combining the estimates above, we get (\ref{eq:oscestU^Qgen}).
\end{proof}

We will also need the following
\begin{lem} \label{lem:lem4}
Let $z_0 \in \overline{\R^{d+1}_+}$, $R_0 \in (0,\infty), \nu \in [32,\infty), r > \frac{R_0}{\nu}$, and $q_0, \sigma \in (1,\infty)$. Suppose $u \in \mathcal{H}^1_{q_0\sigma}(Q_{\nu r}^+(z_0), \mu; \R^n)$ has $\operatorname{supp}u \subset (s-(\delta R_0)^2, s) \times \R^d_+$ for some $s \in \R$ and $\delta >0$. Then
\begin{equation*}
\begin{aligned}
&\fint_{Q_r^+(z_0)}|U| d\mu \leq \left(\fint_{Q_r^+(z_0)}\chi_{(s-(\delta R_0)^2, s)}(t) d\mu\right)^{1-\frac{1}{q_0\sigma}} \left(\fint_{Q_r^+(z_0)}|U|^{q_0\sigma} d\mu\right)^{\frac{1}{q_0\sigma}} \\
&\leq \left(\frac{\delta R_0}{r}\right)^{2-\frac{2}{q_0\sigma}} \left(\fint_{Q_r^+(z_0)}|U|^{{q_0}\sigma} d\mu\right)^{\frac{1}{{q_0}\sigma}} \leq (\delta \nu)^{2-\frac{2}{q_0\sigma}} \left(\fint_{Q_r^+(z_0)}|U|^{{q_0}\sigma} d\mu\right)^{\frac{1}{{q_0}\sigma}}.
\end{aligned}
\end{equation*}
\end{lem}

We use the following filtration of partitions on $\Omega_T$:
\[
\mathbb{C}_m^+ := \{Q^m = Q^m_{(i_0,i_1,\dots,i_d)}: (i_0,i_1,\dots,i_d) \in \mathbb{Z}^{d+1}, i_d \geq 0\},
\]
where $m \in \mathbb{Z}$ and
\[
Q^m_{(i_0,i_1,\dots,i_d)} = (i_0 2^{-2m}, (i_0+1)2^{-2m}) \times (i_0 2^{-m}, (i_0+1)2^{-m}) \times \dots \times (i_d 2^{-m}, (i_d+1)2^{-m})
\]
if $T = \infty$, or
\[
Q^m_{(i_0,i_1,\dots,i_d)} = (T-(i_0+1) 2^{-2m}, T-i_0 2^{-2m}) \times (i_0 2^{-m}, (i_0+1)2^{-m}) \times \dots \times (i_d 2^{-m}, (i_d+1)2^{-m})
\]
if $T < \infty$, where also $i_0 \geq 0$.

We will use the following generalized version of the weighted Fefferman-Stein inequality in the spirit of \cite[Theorem 2.7]{Krylov2009JFA}, which is a special case of \cite[Corollary 2.8]{DongKimTAMS2018}.

\begin{thm} \label{thm:GFS}
Let $\alpha \in (-1,\infty), p,q\in (1,\infty), K,K_4 \in [1,\infty)$, $w = w_0(t)w_1(x)$, where $w_0 \in A_q(\R), w_1 \in A_p(\R^d_+,\mu)$ with $[w_0]_{A_q(\R)} \leq K, [w_1]_{A_p(\R^d_+,\mu)} \leq K$, and $f,g \in L_{q,p}(\Omega_T, wd\mu)$. Suppose that for each $m \in \mathbb{Z}$ and $Q^m \in \mathbb{C}_m^+$ there exists a measurable function $f^C$ on $Q^m$ such that $\frac{1}{K_4}|f| \leq f^C \leq K_4 |f|$ a.e. on $Q^m$, and
\[
\fint_{Q^m}|f^C(z)-(f^C)_{Q^m, \mu}|d\mu(z) \leq g(y) \quad \forall y \in Q^m.
\]
Then
\[
\Vert f \Vert_{L_{q,p}(\Omega_T, wd\mu)} \leq N \Vert g \Vert_{L_{q,p}(\Omega_T, wd\mu)},
\]
where $N=N(d,\alpha,p,q,K,K_4)>0$.
\end{thm}

\begin{proof}[Proof of Theorem {\rm \ref{thm:main1} (i)}]
For the proof, we use the strategy from \cite{DongKimTAMS2018}. We first show the a-priori estimate (\ref{eq:estmain1}). Let $F^i,f \in L_{q,p}(\Omega_T,w d\mu; \R^n)$ and $u \in \mathcal{H}^1_{q,p}(\Omega_T,w d\mu; \R^n)$ satisfy \eqref{eq:eqpar}-\eqref{eq:conormal}. Assume for a moment that $B^i = \Hat{B}^i = C = 0$.

First, from the reverse H\"older property of Muckenhoupt weights, there exists $p_0=p_0(d,\alpha,p,q,K) \in (1,\min(p,q))$ sufficiently close to 1 such that
\begin{equation} \label{eq:Apimprove}
w_0 \in A_{q/p_0}(\R), \quad w_1 \in A_{p/p_0}(\R^d_+, \mu),
\end{equation}
and therefore,
\[
L_{q,p, \text{loc}}(\overline{\R^{d+1}_+},w d\mu) \subset L_{p_0, \text{loc}}(\overline{\R^{d+1}_+},\mu).
\]
Moreover, $[w_1]_{A_{p/p_0}(\R^d_+,\mu)} \leq \Tilde{K} = \Tilde{K}(d,\alpha,p,K)$ and $[w_0]_{A_{q/p_0}(\R)} \leq \Tilde{K} = \Tilde{K}(q,K)$. For the proof, see \cite[Theorem 7.2.5]{GrafakosCFA}. Fix any $q_0>1$ and $\sigma>1$ with $q_0\sigma=p_0$; for instance, take $q_0=\sigma=p_0^{1/2}$.

We first consider the case $0 < R_0 < \infty$. Let $\delta > 0$ and $\nu \in [32,\infty)$ be some constants which will be determined later. For the time being, assume that $\operatorname{supp}u \subset (s-(\delta R_0)^2, s) \times \R^d_+$ for some $s \in \R$. Fix any $Q^m \in \mathbb{C}_m^+$. We claim that there exists $U^C$ on $Q^m$ such that $N^{-1}_{\kappa}U \leq U^C \leq N_{\kappa}U$ on $Q^m$ and for any $z \in Q^m$ it holds that

\begin{equation} \label{eq:est2}
\begin{gathered}
\fint_{Q^m}|U^C-(U^C)_{Q^m}| d\mu \leq N \nu^{\frac{d+2+\alpha_+}{{q_0}}}(\mathcal{M}^{1/{q_0}}(|F|^{q_0})(z)+\mathcal{M}^{1/{q_0}}(|f|^{q_0})(z)) \\
+ N(\nu^{-1}+\nu^{\frac{d+2+\alpha_+}{{q_0}}}\gamma_0^{\frac{1}{{q_0}\sigma'}}+(\delta \nu)^{2-\frac{2}{q_0\sigma}}) \mathcal{M}^{1/({q_0}\sigma)}(|U|^{{q_0}\sigma})(z).
\end{gathered}
\end{equation}
Let $z_0 \in \overline{Q^m}$ be the top-right corner of $\overline{Q^m}$ and let $r=2^{\frac{d-1}{2}}2^{-m}$ so that $Q^m \subset Q_r^+(z_0)$. If $r \in (0,\frac{R_0}{\nu}]$, then (\ref{eq:est2}) follows by Proposition \ref{prop:prop3} (we take $U^C = U^Q$, where $U^Q$ is from Proposition \ref{prop:prop3}). If $r > \frac{R_0}{\nu}$, we take $U^C = U$ and (\ref{eq:est2}) follows by Lemma \ref{lem:lem4}.

Using (\ref{eq:Apimprove}), we apply the generalized weighted Fefferman-Stein inequality (Theorem \ref{thm:GFS}) and the weighted Hardy-Littlewood inequality (Theorem \ref{thm:HLFS}), and get
\begin{equation*}
\begin{gathered}
\Vert U \Vert_{L_{q,p}(\Omega_T,w d\mu)} \leq N \nu^{\frac{d+2+\alpha_+}{{q_0}}}(\Vert F \Vert_{L_{q,p}(\Omega_T,w d\mu)} + \Vert f \Vert_{L_{q,p}(\Omega_T,w d\mu)}) \\
 + N(\nu^{-1}+\nu^{\frac{d+2+\alpha_+}{{q_0}}}\gamma_0^{\frac{1}{{q_0}\sigma'}}+(\delta \nu)^{2-\frac{2}{q_0\sigma}})\Vert U \Vert_{L_{q,p}(\Omega_T,w d\mu)}.
\end{gathered}
\end{equation*}
Fixing $\nu = \max(32, 4N)$, $\gamma_0>0$ with $\nu^{\frac{d+2+\alpha_+}{{q_0}}}\gamma_0^{\frac{1}{{q_0}\sigma'}} = 1/(8N)$, and $\delta>0$ with $(\delta\nu)^{2-\frac{2}{q_0\sigma}}=1/(8N)$, we arrive to (\ref{eq:estmain1}). Note that so far we did not impose any restriction on $\lambda$ other than $\lambda \geq 0$.

To remove the assumption $\operatorname{supp}u \subset (s-(\delta R_0)^2, s) \times \R^d_+$, one can use a standard partition of unity in time argument. Precisely, fix a cut-off function $\xi=\xi(t) \in C_0^{\infty}((-(R_0 \delta)^2, 0))$ with $\int_{\R}\xi(s)^q ds = 1$ and $\int_{\R} |\xi'(s)|^q ds \leq \frac{N}{(R_0 \delta)^{2q}}$. For each $s \in \R_T$, apply estimate (\ref{eq:estmain1}) to $u(t,x)\xi(t-s)$, raise it to power $q$ and integrate over $s \in \R_T$. Here one needs $\lambda \geq \frac{N}{(\delta R_0)^2}$ to close the argument. We set $\lambda_0 := \frac{N}{\delta^2}$. Note that $\delta$, all constants $N$, and $\lambda_0$ depend on $d,\alpha,p,q,\kappa,$ and $K$.

If $R_0 = \infty$, assumption $\operatorname{supp}u \subset (s-(\delta R_0)^2, s) \times \R^d_+$ is not needed and estimate (\ref{eq:est2}) holds without the term $(\delta\nu)^{2-\frac{2}{q_0\sigma}}$. Hence, if $R_0 = \infty$, estimate (\ref{eq:estmain1}) holds for any $\lambda \geq 0$.

In the general case, we treat the lower-order terms with $B,\Hat{B}$, and $C$ as perturbations. Let $N_0$ be the constant in estimate (\ref{eq:estmain1}) for $B^i=\Hat{B}^i=C=0$ and let $\lambda \geq \max(\lambda_0 R_0^{-2}, MK_1^2)$, where the constant $M>0$ will be determined below. Since $|B| \leq K_1x_d, |\Hat{B}| \leq K_1x_d, |C| \leq K_1^2 x_d^2$, the terms $D_i(x_d^{\alpha-1}B^i u)$ and $x_d^{\alpha-1}\Hat{B}^i D_i u+x_d^{\alpha-2}Cu$ can be combined with respectively $D_i(x_d^{\alpha}F^i)$ and $\sqrt{\lambda}x_d^{\alpha} f$ on the right-hand side of the equation. This gives ($\Vert \cdot \Vert = \Vert \cdot \Vert_{L_{q,p}(\Omega_T,w d\mu)}$)
\[
\begin{gathered}
\Vert Du \Vert + \lambda^{1/2} \Vert u \Vert \leq N_0 \Vert F \Vert + N_0 \Vert f \Vert + N_0 K_1 \Vert u \Vert + \lambda^{-1/2} N_0 K_1 \Vert Du \Vert + \lambda^{-1/2} N_0 K_1^2 \Vert u \Vert \\
\leq N_0 \Vert F \Vert + N_0 \Vert f \Vert + N_0 M^{-1/2} \Vert Du \Vert + N_0 K_1 (M^{-1/2}+1) \Vert u \Vert.
\end{gathered}
\]
Taking $M$ sufficiently large, for example $M = 16N_0^2+1$, and $N=2N_0$, gives estimate (\ref{eq:estmain1}) in full generality.

The existence of solutions can be proved by adapting the argument in \cite[Section 8]{DongKimTAMS2018}. We give a sketch of the proof. It suffices to consider $T=\infty$. By the method of continuity, it suffices to consider $A^{ij}=\delta^{ij}I_{n \times n}$ and $C_0 = I_{n \times n}$. We keep $a_0(x_d)$ unchanged as the definition of the space $\mathcal{H}^1_{q,p}$ depends on $a_0(x_d)$. Again from the reverse H\"older property of Muckenhoupt weights, see \cite[Theorem 7.2.5]{GrafakosCFA}, there exists large enough $p_1 \in (\max(p,q),\infty)$ such that
\[
w_0^{\frac{p_1}{p_1-q}}(t) \in A_q(\R), \quad w_1^{\frac{p_1}{p_1-p}}(x) \in A_p(\R^d_+, \mu).
\]
In particular, the measures $w_0^{\frac{p_1}{p_1-q}}(t)dt$ and $w_1^{\frac{p_1}{p_1-p}}(x) d\mu(x)$ are doubling, and
\[
L_{p_1, \text{loc}}(\overline{\R^{d+1}_+}, \mu) \subset L_{q,p,\text{loc}}(\overline{\R^{d+1}_+}, wd\mu).
\]
Take sequences $F_k, f_k \in C_0^{\infty}(\R^{d+1}_+)$ such that $F_k \to F$ and $f_k \to f$ in $L_{q,p}(\R^{d+1}_+, wd\mu)$. By Theorem \ref{thm:exisuniqsimpleunweighted}, there exist solutions $u_k \in \mathcal{H}^1_{p_1}(\R^{d+1}_+, \mu; \R^n)$ with $F_k, f_k$ in place of $F,f$. By using the doubling property of measures $w_0^{\frac{p_1}{p_1-q}}(t)dt$ and $w_1^{\frac{p_1}{p_1-p}}(x) d\mu(x)$, an iteration argument similar to the one used in the proof of Lemma \ref{lem:LptoLsmall} shows that $u_k \in \mathcal{H}^1_{q,p}(\R^{d+1}_+, wd\mu; \R^n)$. See also \cite[Section 8]{DongKimTAMS2018} for details on the iteration argument in the mixed-norm weighted setting. By the a-priori estimate (\ref{eq:estmain1}), the sequence $u_k$ is Cauchy in $\mathcal{H}^1_{q,p}(\R^{d+1}_+, wd\mu; \R^n)$, and therefore converges to the solution $u \in \mathcal{H}^1_{q,p}(\R^{d+1}_+, wd\mu; \R^n)$ with $F,f$.
\end{proof}

\section{Equations with singular lower order terms} \label{sec:lowerterms}

In this section, we finish the proof of Theorem \ref{thm:main1} and prove Theorem \ref{thm:mainYp}. We start with
\begin{proof}[Proof of Remark {\rm \ref{rem:1} (3)}]
We only need to show the $\Rightarrow$ implications.

Suppose that $w \in A_p(\R_+,\mu) \cap X_p(\mu)$. By H\"older's inequality,
\[
\infty > \left(\int_0^1 w(y)y^{-p}d\mu(y)\right)^{1/p}\left(\int_0^1 w(y)^{-\frac{1}{p-1}} d\mu(y)\right)^{1/p'} \geq \int_0^1 y^{-1}d\mu(y) = \int_0^1 y^{\alpha-1}dy.
\]
Hence, $\alpha > 0$. By duality and the above,
\[
\exists w \in A_p(\R_+,\mu) \cap Y_p(\mu) \iff \exists w^{-\frac{1}{p-1}} \in A_{p'}(\R_+,\mu) \cap X_{p'}(\mu) \iff \alpha > 0.
\]

Now suppose that $w \in X_p(\mu) \cap Y_p(\mu)$. By H\"older's inequality,
\[
\infty > \left(\int_0^1 w(y)y^{-p}d\mu(y)\right)^{1/p}\left(\int_0^1 w(y)^{-\frac{1}{p-1}}y^{-p'} d\mu(y)\right)^{1/p'} \geq \int_0^1 y^{-2}d\mu(y) = \int_0^1 y^{\alpha-2}dy.
\]
Hence, $\alpha > 1$.
\end{proof}

In the following two lemmas, we give equivalent characterizations of the $X_p(\mu)$ and $Y_p(\mu)$ conditions in terms of certain weighted Hardy inequalities.
\begin{lem} \label{lem:Hardyinfty}
Let $\alpha \in \R, d\mu(y) = y^{\alpha} dy$ for $y \in \R_+$, $p \in (1, \infty)$, and $w$ a weight on $\R_+$. Then, there exists $N>0$ such that
\[
\left( \int_0^{\infty} \left|y^{-1}\int_y^{\infty} f(s)ds\right|^p w(y) d\mu(y)\right)^{1/p} \leq N \left( \int_0^{\infty} |f(y)|^p w(y)  d\mu(y)\right)^{1/p}
\]
holds for all $f \in L_p(\R_+,wd\mu)$ if and only if $w\in X_p(\mu)$. Moreover, the least constant $N$ in the above inequality satisfies $[w]_{X_p(\mu)} \leq N \leq 2[w]_{X_p(\mu)}$.
\end{lem}

\begin{lem} \label{lem:Hardy0}
Let $\alpha \in \R, d\mu(y) = y^{\alpha} dy$ for $y \in \R_+$, $p \in (1, \infty)$, and $w$ a weight on $\R_+$. Then, there exists $N>0$ such that
\[
\left(\int_0^{\infty}\left|y^{-\alpha}\int_0^y s^{\alpha-1} f(s)ds\right|^p w(y)dy\right)^{1/p} \leq N \left(\int_0^{\infty}|f(y)|^p w(y)dy\right)^{1/p}
\]
holds for all $f \in L_p(\R_+,wd\mu)$ if and only if $w\in Y_p(\mu)$. Moreover, the least constant $N$ in the above inequality satisfies $[w]_{Y_p(\mu)} \leq N \leq 2[w]_{Y_p(\mu)}$.
\end{lem}
Lemma \ref{lem:Hardyinfty} and Lemma \ref{lem:Hardy0} are special cases of \cite[Theorem 2]{Muckenhoupt1972} and \cite[Theorem 1]{Muckenhoupt1972}, respectively.

\begin{proof}[Proof of Theorem {\rm \ref{thm:main1} (ii)}]
It remains to prove the a-priori estimate (\ref{eq:estmainXp}) as the existence of solutions would then follow from the method of continuity and the existence from Part (i).

As $w_3 \in X_p(\mu)$, by Hardy's inequality (Lemma \ref{lem:Hardyinfty}), we have
\begin{equation} \label{eq:HardyXp}
\Vert u/x_d \Vert_{L_{q,p}(\Omega_T,w d\mu)} \leq 2[w_3]_{X_p(\mu)} \Vert D_d u \Vert_{L_{q,p}(\Omega_T,w d\mu)}.
\end{equation}
We now split
\begin{gather*}
D_i(x_d^{\alpha-1}B^i u) = D_i(x_d^{\alpha}(B^i/x_d) I_{|B| \leq K_1x_d} u) + D_i(x_d^{\alpha}B^i I_{K_1x_d < |B| \leq \varepsilon K_2^{-1}} u/x_d) = {\rm I} + {\rm II}, \\
x_d^{\alpha-2}C u = x_d^{\alpha} (C/x_d^2) I_{|C| \leq K_1^2 x_d^2}u + x_d^{\alpha} (C/x_d) I_{K_1^2 x_d^2 < |C| \leq K_1 K_2^{-1} x_d} u/x_d = {\rm I}'+{\rm II}'.
\end{gather*}
Let $N_1>0$ be the constant in estimate (\ref{eq:estmain1}) and let $M_1$ be the constant $M$ from Theorem \ref{thm:main1} (i). Assume that $\lambda \geq \max(\lambda_0R_0^{-2},MK_1^2)$, where the constant $M \geq M_1$ will be determined below. By combining the terms ${\rm II}$ and ${\rm II}'$ with respectively $D_i(x_d^{\alpha}F^i)$ and $\sqrt{\lambda}x_d^{\alpha}f$ in the right hand side of the equation, using (\ref{eq:estmain1}) and (\ref{eq:HardyXp}), we get ($\Vert \cdot \Vert = \Vert \cdot \Vert_{L_{q,p}(\Omega_T,w d\mu)}$)
\begin{equation*}
\begin{gathered}
\Vert Du \Vert + \sqrt{\lambda} \Vert u \Vert \leq N_1 \Vert F \Vert + N_1 \Vert f \Vert + N_1 (\varepsilon K_2^{-1}+\lambda^{-1/2}K_1 K_2^{-1}) \Vert u/x_d \Vert \\
\leq N_1 \Vert F \Vert+ N_1 \Vert f \Vert + 2N_1 (\varepsilon+M^{-1/2}) \Vert D_d u \Vert.
\end{gathered}
\end{equation*}

Taking $\varepsilon = 1/(8N_1)$, $M=\max(64N_1^2,M_1)$, and using (\ref{eq:HardyXp}), we arrive to (\ref{eq:estmainXp}) with $N=6N_1$.
\end{proof}

\begin{proof}[Proof of Theorem {\rm \ref{thm:mainYp}}]
(i) Define $G_0$ by $\lim_{x_d \to 0^+} x_d^{\alpha}G_0 = 0$ and
\[
x_d^{\alpha-1}F_0 = D_d(x_d^{\alpha} G_0).
\]
That is,
\[
G_0(\cdot,x_d) = x_d^{-\alpha}\int_0^{x_d}s^{\alpha-1}F_0(\cdot,s)ds.
\]
As $w_3 \in Y_p(\mu)$, by Lemma \ref{lem:Hardy0} it holds that
\begin{equation} \label{eq:Hardy0}
\Vert G_0 \Vert_{L_{q,p}(\Omega_T,w d\mu)} \leq 2[w_3]_{Y_p(\mu)} \Vert F_0 \Vert_{L_{q,p}(\Omega_T,w d\mu)}.
\end{equation}
Then for any $\varphi \in C_0^{\infty}((-\infty,T)\times \R^d; \R^n)$, we have
\begin{equation} \label{eq:5}
\int_{\Omega_T}(F^i \cdot D_i \varphi + x_d^{-1} F_0 \cdot \varphi + \sqrt{\lambda} f\cdot \varphi) d\mu =  \int_{\Omega_T}(F^i \cdot D_i \varphi - G_0 \cdot D_d \varphi + \sqrt{\lambda} f\cdot \varphi)d\mu.
\end{equation}
Hence, both sides of (\ref{eq:soldefvar}), viewed as functionals on $C_0^{\infty}((-\infty,T)\times \R^d; \R^n)$, are indeed elements of $\mathbb{H}^{-1}_{q,p}(\Omega_T,wd\mu; \R^n)$. This justifies the application of the method of continuity between $\mathcal{H}^1_{q,p}(\Omega_T,wd\mu; \R^n)$ and $\mathbb{H}^{-1}_{q,p}(\Omega_T,wd\mu; \R^n)$. It suffices to prove the a-priori estimate (\ref{eq:estmainYp}) as the existence of solutions would then follow from the method of continuity and the existence from Theorem \ref{thm:main1} (i).

Let $N_1>0$ be the constant in estimate (\ref{eq:estmain1}) and let $M_1$ be the constant $M$ from Theorem \ref{thm:main1} (i). Assume that $\lambda \geq \max(\lambda_0R_0^{-2},MK_1^2)$, where the constant $M \geq M_1$ will be determined below. We split
\begin{gather*}
x_d^{\alpha-1}\Hat{B}^i D_i u = x_d^{\alpha}(\Hat{B}^i/x_d) I_{|\Hat{B}| \leq K_1x_d} D_i u + x_d^{\alpha-1}\Hat{B}^i I_{K_1x_d < |\Hat{B}| \leq \varepsilon K_3^{-1}} D_i u = {\rm I} + {\rm II}, \\
x_d^{\alpha-2}C u = x_d^{\alpha} (C/x_d^2) I_{|C| \leq K_1^2 x_d^2}u + x_d^{\alpha-1} (C/x_d) I_{K_1^2 x_d^2 < |C| \leq K_1 K_3^{-1} x_d} u = {\rm I}'+{\rm II}'.
\end{gather*}
By combining the terms ${\rm II}$ and ${\rm II}'$ with $x_d^{\alpha-1}F_0$ on the right-hand side of the equation, using  (\ref{eq:5}) and estimates (\ref{eq:Hardy0}), (\ref{eq:estmain1}), we have ($\Vert \cdot \Vert = \Vert \cdot \Vert_{L_{q,p}(\Omega_T,w d\mu)}$)
\[
\begin{gathered}
\Vert Du \Vert + \sqrt{\lambda} \Vert u \Vert \leq N_1 \Vert F \Vert + N_1 \Vert f \Vert + 2N_1 K_3(\Vert F_0 \Vert + \varepsilon K_3^{-1}\Vert Du \Vert +K_1 K_3^{-1} \Vert u \Vert).
\end{gathered}
\]

Taking $\varepsilon=1/(4N_1)$, $M=\max(16N_1^2,M_1)$, and $N=4N_1$, we arrive to estimate (\ref{eq:estmainYp}).

(ii) Again, we only need to prove the a-priori estimate (\ref{eq:estmainXpYp}). As $w_3 \in X_p(\mu) \cap Y_p(\mu)$, we can derive estimate (\ref{eq:estmainXpYp}) with appropriate constants $\varepsilon,M$, and $N$, by combining Part (i), Theorem \ref{thm:main1} (ii), and splitting
\[
x_d^{\alpha-2}Cu = x_d^{\alpha-2}C I_{|C|\leq h(x_d)}u + x_d^{\alpha-1}C I_{h(x_d) < |C| \leq \varepsilon(K_2K_3)^{-1}} u/x_d = {\rm I} + {\rm II},
\]
where $h(x_d) := \max(K_1 K_2^{-1}x_d, K_1K_3^{-1}x_d, K_1^2 x_d^2)$, and we absorb the ${\rm II}$ term into $x_d^{\alpha-1}F_0$.
\end{proof}

We finish this Section with the proof of Remark \ref{rem:1} (6). We first present the following lemma that shows that $C_0^{\infty}(\R^d_+)$ is dense in $W^1_p(\R^d_+,wd\mu)$ if $w_3 \in X_p(\mu)$.
\begin{lem} \label{lem:densityXp}
Let $\alpha \in (0,\infty), p \in (1,\infty), w(x)=w_2(x')w_3(x_d)$, where $w_2 \in A_p(\R^{d-1})$ and $w_3 \in A_p(\R_+, \mu) \cap X_p(\mu)$. Let $\zeta \in C_0^{\infty}(\R)$ with $\zeta(s) = 0$ for $s \leq 1$ and $\zeta(s) = 1$ for $s \geq 2$. Then for any $u \in W^1_p(\R^d_+, wd\mu)$,
\[
u_k(x) := u(x)\zeta(kx_d) \to u \quad \text{as } k \to \infty \quad \text{in } W^1_p(\R^d_+, wd\mu).
\]
\end{lem}

\begin{proof}
By the dominated convergence theorem
\[
u_k \to u, \quad D_{x'} u_k \to D_{x'} u \quad \text{in } L_p(\R^d_+, wd\mu).
\]
By Lemma \ref{lem:Hardyinfty}, $u/x_d \in L_p(\R^d_+, wd\mu)$. Then
\begin{gather*}
\Vert D_d u_k - D_d u \Vert_{L_p(\R^d_+, wd\mu)} = k \Vert u \zeta'(kx_d) \Vert_{L_p(\R^d_+, wd\mu)} \leq N k \Vert u I_{x_d \leq 2/k} \Vert_{L_p(\R^d_+, wd\mu)} \\
\leq N \Vert u/x_d I_{x_d \leq 2/k} \Vert_{L_p(\R^d_+, wd\mu)} \to 0
\end{gather*}
as $k \to \infty$ by the dominated convergence theorem.
\end{proof}

We now give the proof of Remark \ref{rem:1} (6). That is, in Theorem \ref{thm:mainYp} from $w_3(x_d) \in Y_p(\mu)$ it follows that any interior weak solution $u$ to (\ref{eq:eqpar}) with $u, Du, F,f \in L_{q,p}(\Omega_T, wd\mu)$  automatically satisfies the conormal boundary condition (\ref{eq:conormal}). Indeed, fix any $\varphi \in C_0^{\infty}((-\infty,T) \times \R^d; \R^n)$. Denote by $\Tilde{w}(t,x)=w_0(t)^{-\frac{1}{q-1}} w_2(x')^{-\frac{1}{p-1}} w_3(x_d)^{-\frac{1}{p-1}}$ the weight dual to $w$. By H\"older's inequality, for any $f,g$, it holds that 
\[
\left|\int_{\Omega_T}f\cdot g d\mu\right| \leq \Vert f \Vert_{L_{q,p}(\Omega_T, wd\mu; \R^n)} \Vert g \Vert_{L_{q',p'}(\Omega_T, \Tilde{w}d\mu; \R^n)}.
\]
As mentioned in Remark \ref{rem:1} (2), $w_3(x_d)^{-\frac{1}{p-1}} \in X_{p'}(\mu)$. We have that (\ref{eq:soldefvar}) holds with $\varphi_k$ instead of $\varphi$, where $\varphi_k=\varphi(t,x)\zeta(kx_d) \in C_0^{\infty}(\Omega_T)$ is the approximating sequence from Lemma \ref{lem:densityXp}. By Lemma \ref{lem:densityXp}, $\partial_t \varphi_k, D\varphi_k, \varphi_k/x_d, \varphi_k$ converge in $L_{q',p'}(\Omega_T, \Tilde{w}d\mu)$ to respectively $\partial_t \varphi, D\varphi, \varphi/x_d, \varphi$. Passing to the limit as $k \to \infty$, identity (\ref{eq:soldefvar}) holds for any $\varphi \in C_0^{\infty}((-\infty,T) \times \R^d; \R^n)$.

\section{Hilbert space-valued solutions} \label{sec:7}

We finish the paper by discussing the extension of our results to the general case of Hilbert space-valued solutions. Let $H$ be a Hilbert space over $\R$ or $\C$, not necessarily separable. Let $\langle \cdot, \cdot \rangle_H$ denote the Hermitian form on $H$, which we assume is antilinear in the second coordinate. Let $\mathcal{B}(H)$ be the algebra of bounded linear operators on $H$. If $a \in \mathcal{B}(H)$, $a^*$ stands for the operator adjoint to $a$.

Throughout the paper, all of our main and intermediate results continue to hold almost verbatim with $\R^n$ changed to $H$, $\R^{n\times n}$ changed to $\mathcal{B}(H)$, and the following modifications in mind. All dot products $u\cdot v$ should be understood as $\langle u,v\rangle_H$. All vector-valued functions $u,D_i u, F^i,F_0, f: \Omega_T \to H$ etc. are strongly measurable in the sense of Bochner and all integrals of vector-valued functions are Bochner integrals. We refer to \cite[Chapter 1]{AnalysisBanSpacesI} for the definitions of Bochner integral and measurability for vector-valued functions.

We now discuss the measurability conditions on the operator-valued coefficients. Here a distinction between $H$ being separable or non-separable needs to be made. If $H$ is non-separable, we assume that $a_0,A^{ij},B^i,\Hat{B}^i,C,C_0:\Omega_T \to \mathcal{B}(H)$ etc. are strongly measurable with respect to the norm topology on $\mathcal{B}(H)$ (i.e., in the sense of Bochner), which is a usual assumption in the theory of parabolic PDEs in Banach spaces. If $H$ is separable, we impose a weaker assumption that $a_0,A^{ij},B^i,\Hat{B}^i,C,C_0:\Omega_T \to \mathcal{B}(H)$ etc. are strongly measurable with respect to the strong operator topology. That is, for every $v \in H$ the $H$-valued function $C_0v:\Omega_T \to H$ is strongly measurable, and similarly for the other coefficients.  We impose these conditions in order for the scalar functions $|a_0|, |A|, |B|, |\Hat{B}|, |C|, |C_0|$ to be measurable on $\Omega_T$. If $\operatorname{dim}H<\infty$, then $\operatorname{dim}\mathcal{B}(H)<\infty$ and therefore both conditions are equivalent to the usual Lebesgue measurability.

Next, all matrix transposes such as $a_0^T$ are understood as $a_0^*$. If the base field is $\C$, we modify (\ref{eq:ellipticity}) and (\ref{eq:nondega0c0}) as follows:
\[
\kappa\sum_{i=1}^d |\xi|_H^2 \leq \operatorname{Re} \sum_{i,j=1}^d \langle A^{ij}\xi_j, \xi_i \rangle_H, \quad \kappa\sum_{i=1}^d |\xi|_H^2 \leq \operatorname{Re} \langle C_0\xi , \xi \rangle_H.
\]
In the proofs of Theorems \ref{thm:L2solv} and \ref{thm:1/2timeL2} in Section \ref{sec:L2}, we change all integrals to their real parts. Starting from Section \ref{sec:4}, the proofs carry forward verbatim by changing $\R^n$ to $H$ and with the above conventions in mind.

\providecommand{\bysame}{\leavevmode\hbox to3em{\hrulefill}\thinspace}
\providecommand{\MR}{\relax\ifhmode\unskip\space\fi MR }
\providecommand{\MRhref}[2]{%
  \href{http://www.ams.org/mathscinet-getitem?mr=#1}{#2}
}
\providecommand{\href}[2]{#2}

\end{document}